\documentclass[12pt]{article}
\usepackage{amsmath}
\usepackage{mathrsfs}
\usepackage{bbm}
\usepackage{amssymb,a4}
\usepackage{amsfonts,amssymb,amsthm}
\usepackage{extarrows}
\usepackage{color}
\usepackage[all,pdf]{xy}

\allowdisplaybreaks  

\newcommand{\al}{\alpha}
\newcommand{\be}{\beta}
\newcommand{\dt}{\delta}
\newcommand{\lmd}{\lambda}

\newcommand{\p}{\partial}
\newcommand{\sgm}{\sigma}

\newcommand{\ov}{\overline}
\newcommand{\spanc}[1]{\mathrm{span}_\C\{#1\}}
\newcommand{\pfend}{\hfill{}$\Box$\end{pf}}
\newcommand{\f}[2]{{F_{#1,#2}}}
\newcommand{\ve}[1]{{v_{#1}}}
\newcommand{\V}{V(\al,\be,F)}

\newcommand{\Z}{{\mathbb Z}}
\newcommand{\C}{{\mathbb C}}
\newcommand{\zp}{{\Z_+}}
\newcommand{\pz}{p\Z}

\newcommand{\g}{\mathfrak{g}}
\newcommand{\U}{{\mathcal U}}
\newcommand{\gz}{{\g{(0)}}}
\newcommand{\bei}{{\be_i}}
\newcommand{\beis}{{\be_{i+s}}}
\newcommand{\beims}{{\be_{i-s}}}
\newcommand{\beir}{{\be_{i-r}}}

\newcommand{\dercq}{\text{Der}\C_Q}
\newcommand{\bfm}{{\bf m}}
\newcommand{\bfn}{{\bf n}}
\newcommand{\bfr}{{\bf r}}
\newcommand{\bfs}{{\bf s}}

\newcommand{\dif}{\frac{d}{dt}}
\newcommand{\A}{\mathcal{A}}
\newcommand{\omg}[2]{{\Omega_{(#1,#2)}^l}}

\newcommand{\vlmd}{v_\lmd}

\newtheorem{thm}{Theorem}[section]
\newtheorem{prop}[thm]{Proposition}
\newtheorem{lem}[thm]{Lemma}

\newtheorem{definition}[thm]{Definition}

\theoremstyle{remark}
\newtheorem*{pf}{{\bf Proof}}

\begin{document}

\begin{center}
{\Large {\bf Classification of irreducible Harish-Chandra modules over gap-$p$ Virasoro algebras}}\\
\vspace{0.5cm}
\end{center}

\begin{center}
{Chengkang Xu\footnote{
The author is supported by the National Natural Science Foundation of China(No.11626157, 11801375),
the Science and Technology Foundation of Education Department of
Jiangxi Province(No. GJJ161044).}\\
School of Mathematical Sciences, Shangrao Normal University, Shangrao, Jiangxi,
China}
\end{center}

\begin{abstract}
We prove that any irreducible Harish-Chandra module for a class of Lie algebras,
which we call gap-$p$ Virasoro algebras,
must be a highest weight module, a lowest weight module,
or a module of intermediate series.
These algebras are closely related to
the Heisenberg-Virasoro algebra and the algebra of derivations over a quantum torus.
They also contain subalgebras which are isomorphic to the Virasoro algebra $Vir$,
but graded by $\pz$ (unlike $Vir$ by $\Z$).
\end{abstract}

\section{Introduction}

\def\theequation{1.\arabic{equation}}
\setcounter{equation}{0}

Throughout this paper, $\C,\Z,\zp$ refer to the set of complex numbers,
integers, and positive integers respectively.
Let $p>1$ be a positive integer and denote by $\g'$ the Lie algebra with a basis
$\{L_m\mid m\in\Z \}$ and Lie bracket
$$[L_m,L_n]=(n-m)L_{m+n},\ \ \ [L_m,L_r]=rL_{m+r},\ \ \ [L_r,L_s]=0,$$
where $m,n\in \pz$ and $r,s\notin \pz$.
One can easily show that the universal central extension $\g$ of $\g'$
has a $p$-dimensional center $\spanc{C_i\mid 0\leq i\leq p-1}$, and Lie brackets
$$\begin{aligned}
&[L_m,L_n]=(n-m)L_{m+n}+\dt_{m+n,0}\frac{1}{12}\left((\frac{m}{p})^3-(\frac{m}{p})\right)C_0;\\
&[L_m,L_r]=rL_{m+r};\ \ \ [L_r,L_s]=\dt_{r+s,0}rC_{\overline r},
\end{aligned}$$
where $m,n\in \pz, r,s\notin \pz$ and ${\overline r}$ represents the residue of $r$ by $p$.

One reason for the Lie algebra $\g'$ (or $\g$) to be interesting is
its relation with the Heisenberg-Virasoro algebra $HV$,
which was first introduced in \cite{ACKP} and
has a basis $\{x_m,I_r,K_x,K_{xI},K_I\mid m,r\in\Z \}$ subjecting to
\begin{equation}\label{eq1.1}
\begin{aligned}
&[x_m,x_n]=(n-m)x_{m+n}+\dt_{m+n,0}\frac{m^3-m}{12}K_x,\\
&[x_m,I_r]=rI_{m+r}+\dt_{m+r,0}(m^2+m)K_{xI},\ \ [I_r,I_s]=r\dt_{r+s,0}K_I,\\
&[K_x,HV]=[K_{xI},HV]=[K_I,HV]=0.
\end{aligned}
\end{equation}
Clearly $Vir=\spanc{x_m,K_x\mid m\in\Z}$ is a Virasoro algebra,
and one can easily see that $\spanc{x_{m},I_r\mid m\in \pz,r\in\Z\backslash \pz}$
forms a subalgebra of the centerless Heisenberg-Virasoro algebra
which is isomorphic to $\g'$.
Through this isomorphism, the algebra $\g$ can be realized as
the universal central extension of
the algebra of some differential operators of order at most one on $\C[t^{\pm1}]$,
$$\spanc{t^{m+1}\dif\mid m\in \pz}\oplus\spanc{t^r\mid r\in\Z\backslash \pz}.$$
The part consisting operators of order one
generates a subalgebra $\g(0)=\spanc{L_m,C_0\mid m\in\pz}$ of $\g$
that is isomorphic to the Virasoro algebra $Vir$.
The isomorphism is given by
$$x_i\mapsto \frac{1}{p}L_{pi};\ \ \ K_x\mapsto \frac{C_0}{p^2}.$$
Since $\g(0)$ is graded by $\pz$,
we call the Lie algebra $\g$ a {\em gap-$p$ Virasoro algebra}.
An intriguing point is that $\g$ can not be imbedded into $HV$.

The second reason why we take interest in the algebra $\g'$ is that
it can be imbedded into the algebra of derivations over a rational quantum torus.
Let $\C_Q=\C[t_1^{\pm1},t_2^{\pm1},\cdots,t_d^{\pm1}]$ denote the quantum torus
with respect to the $d\times d$ matrix $Q=(q_{ij})$,
satisfying commuting relation $t_it_j=q_{ij}t_jt_i$,
where $1<d\in\zp$ and $q_{ij}$'s are all roots of unity such that
$q_{ii}=1$ and $q_{ij}q_{ji}=1$.
For $\bfn=(n_1,\cdots,n_d)\in\Z^d$ we write $t^\bfn=t_1^{n_1}t_2^{n_2}\cdots t_d^{n_d}$.
Let $\p_i$ denote the degree derivation corresponding to $t_i$.
We recall from \cite{BGK} the map
$\sgm(\bfm,\bfn)=\prod_{1\leq i<j\leq d}q_{ji}^{m_jn_i}$ for $\bfm,\bfn\in\Z^d$
and the set
$$R=\{\bfm\in\Z^d\mid \sgm(\bfm,\bfn)=\sgm(\bfn,\bfm)\text{ for any }\bfn\in\Z^d \}.$$
Then the algebra $\dercq$ of derivations over $\C_Q$ has a basis
$$\{t^\bfm\p_i, t^\bfs\mid 1\leq i\leq d, \bfm\in R,\bfs\in\Z^d\backslash R\}$$
subjecting to the Lie brackets
$$\begin{aligned}
&[t^\bfm\p_i,t^\bfn\p_j]=t^{\bfm+\bfn}(n_i\p_j-n_j\p_i);\\
&[t^\bfm\p_i,t^\bfs]=\sgm(\bfm,\bfs)s_it^{\bfm+\bfs};\\
&[t^\bfr,t^\bfs]=(\sgm(\bfr,\bfs)-\sgm(\bfs,\bfr))t^{\bfr+\bfs},
\end{aligned}$$
where $1\leq i,j\leq d,\ \bfm,\bfn\in R$,\ $\bfr,\bfs\in\Z^d\backslash R$.
Choose $\bfs\in\Z^d\backslash R$
whose image in the quotient group $\Z^d/R$ has order $p$ (assume such $\bfs$ exists).
Then the algebra $\g'$ is isomorphic to the subalgebra of $\dercq$ spanned by
$$\{t^{kp\bfs}\p_1,t^{l\bfs}\mid k\in\Z,l\in\Z\backslash \pz\},$$
through the map defined by
$$L_{kp}\mapsto\frac{1}{s_1}t^{kp\bfs}\p_1,\ \ L_l\mapsto t^{l\bfs}
\text{ for }k\in\Z, l\notin\pz.$$

In this paper our main concern is the irreducible Harish-Chandra modules over the algebra $\g$.
The classification problem of irreducible Harish-Chandra modules over Lie algebras
is a priority in the representation theory,
and was solved for many infinite dimensional Lie algebras,
such as the Virasoro algebra \cite{M},
the higher rank Virasoro algebra \cite{LZ1},
the Heisenberg-Virasoro algebra \cite{LZ2},
the algebra of derivations over a commuting torus \cite{BF},
some Lie algebras of Block type \cite{WT,GGS},
and so on.
Like the algebra $\g$,
many of these algebras contains a subalgebra isomorphic to the Virasoro algebra.
The difference is that these subalgebras are all graded by $\Z$,
while the subalgebra $\g(0)$ of $\g$ is graded by $\pz$.
This causes a major trouble for the classification of irreducible Harish-Chandra modules over $\g$
for the lack of Virasoro elements in the grading spaces $\g_s, s\notin\pz$.

The paper is organized as follows.
In Section 2 we study irreducible Harish-Chandra modules for $\g$,
and prove that any such module is a highest weight module, a lowest weight module,
or a cuspidal module.
In Section 3 we further prove that any irreducible cuspidal module $M$ for $\g$
must be a module of intermediate series.
The method used here was originally introduced in \cite{BF},
and we modify it to suit for our setting.
The critical notion is the $\mathcal{A}$-cover $\hat{M}$ of $M$.
It turns out that $\hat{M}$ is also cuspidal and $M$ is a subquotient of $\hat{M}$.
Then we prove that an irreducible cuspidal $\g$-module with associative $\mathcal{A}$-action
has weight multiplicities no more than one.
Section 4 is devoted to the irreducible modules of intermediate series over $\g$,
which we prove must be of the form
$\V$ with parameters $\al,\be\in\C$ and $(p-1)\times p$ matrix $F$ (see Theorem \ref{thm4.3}).
Here, the modules of intermediate series over $\gz$ play an important role,
and severe computations are involved using a linkage method.
In the last section some examples of module of intermediate series with small $p$ are given.

\section{Harish-Chandra modules over the algebra $\g$}
\def\theequation{2.\arabic{equation}}
\setcounter{equation}{0}

In this section we study the Harish-Chandra modules for the algebra $\g$
and prove they are either highest weight modules, lowest weight modules or cuspidal modules.

The algebra $\g$ has a $\Z$-gradation $\g=\bigoplus\limits_{i\in\Z}\g_i$,
where
$$\g_i=\begin{cases}
\spanc{ L_i}, &\text{ if } i\neq 0\\
\spanc{L_0,C_j\mid 0\leq j<p}, &\text{ if } i=0.
\end{cases}$$
The subspace $\g_0$ is a Cartan subalgebra of $\g$,
and $\g$ has a triangular decomposition $\g=\g_-\oplus\g_0\oplus\g_+$ where
$\g_\pm=\spanc{L_m\mid\pm m>0}$.
A $\g$-module $V$ is called a {\em weight module} if $\g_0$ acts diagonalizably on $V$.
For any weight module $V$ we have the weight space decomposition
$V=\bigoplus_{\lmd\in\g_0^*}V_\lmd,$
where $\g_0^*=\text{Hom}_\C(\g_0,\C)$ and
$$V_\lmd=\{v\in V\mid av=\lmd(a)v\text{ for all }a\in\g_0\}.$$
The function $\lmd$ is called a {\em weight} provided $V_\lmd\neq 0$ and
the space $V_\lmd$ is called the {\em weight space} corresponding to $\lmd$.
A weight $\g$-module $V$ is called a {\em Harish-Chandra} module if
$\dim V_\lmd< \infty$ for all weights $\lmd$,
called {\em cuspidal}
if there is some $k\in\zp$ such that $\dim V_\lmd<k$ for all weights $\lmd$,
and furthermore called a {\em module of intermediate series} (abbreviate MOIS)
if $\dim V_\lmd\leq1$ for all weight $\lmd$.

The first class of irreducible Harish-Chandra modules for $\g$
are the highest weight modules.
Let $\lmd\in\g_0^*$ and $\C v_\lmd$ be a one dimensional $\g_0$-module defined by
$$av_\lmd=\lmd(a)v_\lmd,\text{ for any }a\in\g_0.$$
Set $\g_+\vlmd=0$ making $\C\vlmd$ a $(\g_0+\g_+)$-module,
and then we have the induced $\g$-module
$M(\lmd)=\U(\g)\otimes_{\U(\g_0+\g_+)}\C\vlmd=\U(\g_-)\vlmd$,
which is called the Verma module and
by the Poincar$\acute{\text{e}}$-Birkhoff-Witt theorem has a basis
$$L_{-m_1}^{p_1}\cdots L_{-m_r}^{p_r}\vlmd,\ \ \ r\geq0,\
 0<m_1<\cdots<m_r\in\Z,p_1,\cdots,p_r\in\zp.$$

Denote $I(\lmd)=\{1\leq i\leq p-1\mid\lmd(C_i)=0\}$.
\begin{prop}\label{prop2.1}
The $\g$-module $M(\lmd)$ is irreducible if and only if $I(\lmd)=\emptyset$.
\end{prop}
\begin{pf}
Suppose $I(\lmd)\neq\emptyset$ and let $i\in I(\lmd)$.
Set $\g_{(i)}=\bigoplus\limits_{k\in\Z}\g_{kp+i}$.
It is easy to see that $\U(\g)\g_{(i)}\vlmd$ is a proper nonzero $\g$-submodule of $M(\lmd)$.

Conversely, suppose $\lmd(C_i)\neq 0$ for all $0<i<p$.
Any nonzero vector in $w\in M(\lmd)$ may be written as a finite sum
\begin{equation}\label{eq2.1}
w=\sum A_{m_1,\cdots,m_r}L_{-m_1}\cdots L_{-m_r}\vlmd,
\end{equation}
where $r\geq 0,\ m_1\geq\cdots\geq m_r\in\zp,\ A_{m_1,\cdots,m_r}\in\C\backslash\{0\}$.

We use induction on the number $\Gamma(w)$ of elements $L_{-s}$
appearing in the expression in equation (\ref{eq2.1}) such that $s\in \pz\cap\zp$.
If $\Gamma(w)=0$, let $F(w)$ be the set of integers $s\in\zp$ such that
$L_{-s}$ appears in equation (\ref{eq2.1}),
and $n_1$ be the largest one in $F(w)$.
Then we have $L_{n_1}w\neq 0$.
Let $n_2$ be the largest integer in the set $F(L_{n_1}w)$ and we have
$L_{n_2}L_{n_1}w\neq 0$.
Repeating this procedure we get $L_{n_k}\cdots L_{n_2}L_{n_1}w\in\C\vlmd\backslash\{0\}$
for some ${n_k}\leq\cdots\leq {n_2}\leq{n_1}$.
Hence $w$ generates $M(\lmd)$.

Now suppose $\Gamma(w)>0$ and any nonzero vector $w'\in M(\lmd)$ with $\Gamma(w')<\Gamma(w)$
generates $M(\lmd)$.
Apply $L_1$ to $w$ and for each summand in equation (\ref{eq2.1}) we have
$$\begin{aligned}
 L_1L_{-m_1}\cdots L_{-m_r}\vlmd&=[L_1,L_{-m_1}\cdots L_{-m_r}]\vlmd\\
  &=\sum_{i=1}^{r} L_{-m_1}\cdots L_{-m_{i-1}}[L_1,L_{-m_i}] L_{-m_{i+1}}\cdots L_{-m_r}\vlmd,
\end{aligned}$$
from which we see that $L_1w\neq0$ and $\Gamma(L_1w)<\Gamma(w)$.
So by the inductional hypothesis $M(\lmd)$ is an irreducible $\g$-module.
\pfend

From the proof of Proposition \ref{prop2.1} we see that
$M(\lmd)$, if reducible, contains a unique maximal proper $\g$-submodule
$$M^\lmd=\begin{cases}
\bigoplus\limits_{i\in I(\lmd)}\U(\g)\g_{(i)}\vlmd &\text{ if }\lmd\neq0
                                                      \text{ and }I(\lmd)\neq\emptyset;\\
\U(\g)\vlmd\backslash\C\vlmd                &\text{ if }\lmd =0.
\end{cases}$$
So $L(\lmd)=M(\lmd)/M^\lmd$ is an irreducible Harish-Chandra $\g$-module,
called irreducible highest weight module with respect to the highest weight $\lmd$.

We may construct the irreducible lowest weight $\g$-modules similarly and
get same results as the highest weight modules.

Finally we have the following
\begin{thm}\label{thm2.2}
An irreducible Harish-Chandra $\g$-module is a highest weight module,
a lowest weight module, or a cuspidal module.
\end{thm}
\begin{pf}
Suppose that $V$ is an irreducible Harish-Chandra $\g$-module,
and $V$ has no highest or lowest weight.
Let $V=\bigoplus\limits_{i\in\Z}V_{\al+i}$ for some $\al\in\g_0^*$
be the weight space decomposition.
For any integer $n>0$, consider the subalgebra $\mathfrak{a}$ of $\g$ generated by $2p$ elements
$$L_n,L_{n+1},\cdots,L_{n+2p-1}.$$
Clearly $\mathfrak{a}$ has a finite codimension in $\g_+$,
which means that there exists some positive integer $N$ satisfying
$$\g_{\geq N}=\spanc{L_k\mid k\geq N}\subseteq\mathfrak{a}.$$
We claim that
$$\mathcal{K}=\bigcap_{i=0}^{2p-1}\ker L_{n+i}\mid_V=0.$$
Indeed, suppose otherwise and let $0\neq v\in\mathcal{K}$,
then $\g_{\geq N}v\subseteq\mathfrak{a}v=\{0\}$.
Lemma 1.6 in \cite{M} implies that
$V$ has a highest weight, which is a contradiction.
From the claim we have
$$\dim V_{-n+\al}\leq \sum_{i=0}^{2p-1}\dim V_{i+\al}.$$
Similarly we get $\dim V_{n+\al}\leq \sum\limits_{i=0}^{2p-1}\dim V_{-i+\al}$.
\pfend

\section{Classification of irreducible Harish-Chandra modules}
\def\theequation{3.\arabic{equation}}
\setcounter{equation}{0}

In this section we prove
\begin{thm}\label{thm3.1}
An irreducible Harish-Chandra $\g$-module is a highest weight module,
a lowest weight module, or a module of intermediate series.
\end{thm}
To do this we only need to show that any irreducible cuspidal $\g$-module
must have weight multiplicity no more than one by Theorem \ref{thm2.2}.
From now on we fix a cuspidal module $M$ over $\g$ (not necessarily irreducible).
We consider $M$ as a cuspidal module over $\gz$,
which is isomorphic to the Virasoro algebra $Vir$.
Then by Proposition (II.7) in \cite{MP} we see that $C_0M=0$.

\begin{lem}[\cite{BL}]\label{lem3.2}
The dimensions of the homogeneous summands of a
nontrivial $\Z$-graded module for an infinite dimensional Heisenberg algebra on which a
central element acts as nonzero scalar are unbounded.
\end{lem}

For $1\leq i\leq p-1$, denote
$$\g(i)=\spanc{L_{i+pk},C_i\mid k\in\Z}.$$
This is an infinite dimensional Heisenberg algebra,
and $M$ may be considered as a cuspidal $\g(i)$-module.
Then Lemma \ref{lem3.2} implies $C_iM=0$.
So $M$ reduces to a $\g'$-module.

Recall that the algebra $\g'$ has a realization as a part of differential operators
of order at most one on $\C[t^{\pm1}]$,
$$\g'\cong \spanc{t^{m+1}\dif, t^s\mid m\in \pz, s\notin \pz}.$$
Denote $\A=\spanc{t^m\mid m\in \pz}$, which is a unital associative algebra with multiplication
$t^mt^n=t^{m+n}$ for $m,n\in \pz$.

\begin{definition}
An {\em $\A\g'$-module} $V$ is a $\g'$-module with a compatible associative $\A$-action,
by which we mean that
$$t^mt^nv=t^{m+n}v,\ \ L_mt^nv-t^nL_mv=nt^{m+n}v,\ \ \ \ t^nL_rv=L_rt^nv$$
for any $m,n\in \pz, r\notin \pz$ and $v\in V$.
\end{definition}

Let $\g'(0)=\spanc{L_m\mid m\in \pz}$.
We note that an $\A\g'$-module $V$ is equivalent to a $(\g'(0)\ltimes\C[t^{\pm1}])$-module
with an associative $\A$-action, if we define
$$t^sv=L_sv\text{ for all }s\notin \pz, v\in V.$$
Denote $\g''=\spanc{L_s\mid s\notin \pz}$, an abelian ideal of $\g'$.
We can define an $\A\g'$-module structure on $\g''$ by
$$L_m\cdot L_s=[L_m,L_s],\ \ \ \ \ t^n\cdot L_s=L_{s+n}$$
for any $m\in\Z, n\in \pz$ and $s\notin \pz$.

\begin{lem}\label{lem3.4}
The tensor space $\g''\otimes M$ is an $\A\g'$-module defined by
$$t^n(L_s\otimes w)=L_{s+n}\otimes w,\ \ \
  a(L_s\otimes w)=[a,L_s]\otimes w+x\otimes a w$$
for $n\in \pz, s\notin \pz, w\in M$ and $a\in \g'$.
\end{lem}
\begin{pf}
We only need to check the compatibility of the $\g'$-module structure
and the $\A$-action on $\g''\otimes M$.
Let $m,n\in \pz$, $r,s\notin \pz$ and $w\in M$. We have
$$\begin{aligned}
 &L_mt^n(L_s\otimes w)-t^n L_m(L_s\otimes w)\\
=&L_m(L_{s+n}\otimes w)-t^n(sL_{s+m}\otimes w+L_s\otimes L_mw)\\
=&(s+n)L_{s+m+n}\otimes w+L_{s+n}\otimes L_mw-(sL_{s+m+n}\otimes w+L_{s+n}\otimes L_mw)\\
=&nL_{s+m+n}\otimes w=nt^{m+n}(L_s\otimes w),
\end{aligned}$$
and
$$\begin{aligned}
 &L_rt^n(L_s\otimes w)-t^n L_r(L_s\otimes w)\\
=&L_r(L_{s+n}\otimes w)-t^n([L_r,L_s]\otimes w+L_s\otimes L_rw)\\
=&([L_r,L_{s+n}]\otimes w+L_{s+n}\otimes L_rw)-L_{s+n}\otimes L_rw=0,
\end{aligned}$$
proving the lemma.
\pfend

Define a map $\pi:\g''\otimes M\rightarrow M$ by $x\otimes w\mapsto xw$,
and set
$$J=\{\sum_s L_s\otimes w_s\mid s\notin \pz, w_s\in M,
        \sum_s L_{s+n}w_s=0\text{ for any }n\in \pz\}.$$
Clearly, $\pi$ is a $\g'$-module homomorphism and $J\subseteq\ker\pi$.

\begin{lem}\label{lem3.5}
The space $J$ is an $\A\g'$-submodule of $\g''\otimes M$.
\end{lem}
\begin{pf}
Let $m\in \pz$. Notice that
for $\sum_sL_s\otimes w_s\in J$,
$$L_m(\sum_sL_s\otimes w_s)=\sum_ssL_{s+m}\otimes w_s+\sum_sL_s\otimes L_mw_s,$$
and, for any $n\in \pz$,
$$\begin{aligned}
&\sum_ssL_{s+m+n}w_s+\sum_sL_{s+n}L_mw_s
   =\sum_ssL_{s+m+n}w_s+\sum_s[L_{s+n},L_m]w_s+\sum_sL_mL_{s+n}w_s\\
&=-n\sum_s L_{s+m+n}w_s+L_m\left(\sum_sL_{s+n}w_s\right)=0.
\end{aligned}$$
This proves $L_m(\sum_sL_s\otimes w_s)\in J$.
Similarly $L_r(\sum_sL_s\otimes w_s)\in J$ for $r\notin \pz$.
So $J$ is a $\g'$-module.
The compatibility with the $\A$-module structure
may be checked the same way as in Lemma \ref{lem3.4}.
\pfend

We call the quotient module $(\g''\otimes M)/J$ the {\em $\A$-cover of $M$},
denoted by $\hat M$.
Moreover, we denote the image of $a\otimes w\in\g''\otimes M$ in $(\g''\otimes M)/J$
by $\psi(a,w)$.

\begin{prop}\label{prop3.6}
The $\A$-cover $\hat M$ of $M$ is also a cuspidal $\A\g'$-module.
\end{prop}
\begin{pf}
The corollary 3.4 in \cite{BF} shows that there exists $l\in\zp$ such that the operators
$$\omg m n=\sum_{i=0}^l(-1)^i\binom li L_{m-in}L_{in},$$
annihilates $M$ for any $m,n\in \pz$.
We fix such an $l$ and set $\Lambda=\{1,2,\cdots,p\}$.
For any weight $\lmd$ of $\hat M$, we have
$$\hat M_\lmd=\spanc{\psi(L_{s+n},v)\mid s\in\Lambda, n\in \pz,v\in M_{\lmd-s-n}}.$$
Let $S$ denote the subspace of $\hat M_\lmd$ spanned by
$$\{\psi(L_{s+n},v)\mid s\in\Lambda, v\in M_{\lmd-s-n}, n\in \pz, |n|\leq lp\},$$
plus $\psi(L_{s_0+n_0},w)$ if $\lmd=s_0+n_0$ for some $s_0\in\Lambda, n_0\in \pz$
and $w\in M_0$.
Clearly $S$ is finite dimensional since $M$ is cuspidal.\\
{\bf Claim}: $\hat M_\lmd=S$.\\
To prove this claim we have to show
$\psi(L_{s+n},v)\in S$ for any $s\in\Lambda, n\in \pz$ and $v\in M_{\lmd-s-n}$.
We use induction on $|n|$.
If $|n|\leq lp$ then the claim is trivial.
Suppose $|n|>lp$. Without lose of generality we assume $n>lp$ (the case $n<-lp$ is similar).
Then the positive numbers $n-p,n-2p,\cdots,n-lp$ are all smaller than $n$.

Now we may assume $\lmd-s-n\neq0$.
Since $L_0v=(\lmd-s-n)v$, we can write $v=L_0w$ for some $w\in M_{\lmd-s-n}$.
Then we have
$$\begin{aligned}
0&=\frac 1s \omg np L_sw=\frac 1s\sum_{i=0}^l(-1)^i\binom li [L_{n-ip}L_{ip},L_s]w\\
 &=\sum_{i=0}^l(-1)^i\binom li L_{n-ip+s}L_{ip}w+\sum_{i=0}^l(-1)^i\binom li L_{n-ip}L_{ip+s}w\\
 &=\sum_{i=0}^l(-1)^i\binom li L_{n-ip+s}L_{ip}w+\sum_{i=0}^l(-1)^i\binom li L_{ip+s}L_{n-ip}w
   +\sum_{i=0}^l(-1)^i\binom li (ip+s)L_{n+s}w\\
 &=\sum_{i=0}^l(-1)^i\binom li L_{n-ip+s}L_{ip}w+\sum_{i=0}^l(-1)^i\binom li L_{ip+s}L_{n-ip}w.
\end{aligned}$$
Here we have used the identity
$$\sum_{i=0}^l(-1)^i\binom li=\sum_{i=0}^l(-1)^i\binom lii=0.$$
So we get
$$L_{n+s}v=L_{n+s}L_0w
  =-\sum_{i=1}^l(-1)^i\binom li L_{n-ip+s}L_{ip}w-\sum_{i=0}^l(-1)^i\binom li L_{ip+s}L_{n-ip}w,$$
which leads to
\begin{equation}\label{eq3.1}
\psi(L_{n+s},v)=-\sum_{i=1}^l(-1)^i\binom li \psi(L_{n-ip+s},L_{ip}w)
  -\sum_{i=0}^l(-1)^i\binom li \psi(L_{ip+s},L_{n-ip}w).
\end{equation}
Notice that $L_{ip}w\in M_{\lmd-s-(n-ip)}, L_{n-ip}w\in M_{\lmd-s-ip}$ and
$|n-ip|<n, ip\leq lp$ for any $0\leq i\leq l$.
We see that the right hand side of (\ref{eq3.1}) lies in $S$ by the inductional assumption,
which proves the claim, and hence the lemma.
\pfend

Now we need to classify all irreducible cuspidal $\A\g'$-modules.

\begin{thm}\label{thm3.7}
Any irreducible cuspidal $\A\g'$-module must be a MOIS over $\g$.
\end{thm}
\begin{pf}
Let $V$ be such a module.
For $0\leq i\leq p-1$ set $V_{(i)}=\bigoplus_{k\in\Z}V_{i+pk}$
(here $i+pk$ is not necessarily the weight of $V_{i+pk}$),
and then $V=\bigoplus_{i=0}^{p-1}V_{(i)}$ (some $V_{(i)}$ might be 0).
We consider $V$ as an $\A\g'(0)$-module.
For any $V_{(i)}\neq 0$, since the $\A$-action on $V$ is associative,
each $t^m, m\in\pz$, is a bijection between the weight spaces $V_{i+pk}, k\in\Z$.
Hence these weight spaces have a same dimension.
Moreover, since all $V_{(i)}$'s are linked by $L_s, s\notin\pz$,
we see that all weight spaces must have a same dimension.
Hence we may write (as isomorphic vector spaces)
$$V=V_0\otimes\left(\bigoplus_{i\in I}t^i\A\right),$$
where $I$ is the subset of $\{0,1,\cdots,p-1\}$ consisting of $i$ such that $V_{(i)}\neq 0$.
Notice that $\g'(0)$ is isomorphic to the centerless Virasoro algebra.
By Theorem 1 from \cite{B} we see that
$V_0$ is an irreducible finite dimensional module for the Lie algebra
$$\g'(0)_+=\spanc{L_m\mid m\geq0, m\in\pz}.$$
For any positive $ m\in\pz$, set $\g'(0)_{\geq m}=\spanc{L_n\mid n\geq m, n\in\pz}$.
Then we have $\g'(0)_{\geq m}V_0=0$ by representation theory of the Virasoro algebra.
Since the quotient $\g'(0)_+/\g'(0)_{\geq m}$ is a finite dimensional solvable Lie algebra,
we see $\dim V_0=1$ by the Lie's Theorem and the irreducibility of $V_0$.
This proves that $V$ is a MOIS over $\g'$.
\pfend

\noindent
{\bf Proof of Theorem \ref{thm3.1}}:
Let $M$ be an irreducible cuspidal $\g$-module.
Consider the decomposition series of the $\A$-cover $\hat M$ as $\A\g'$-module
$$0=\hat M_0\subset\hat M_1\subset\hat M_2\subset\cdots\subset\hat M_l=\hat M$$
where $l\in\zp$ and $\hat M_i/\hat M_{i-1}$ are irreducible $\A\g'$-modules.
Recall the $\g'$-homomorphism $\pi:\hat M\rightarrow M$.
Let $k$ be the largest integer such that $\pi(\hat M_{k-1})=0$.
By the irreducibility of $M$, we see that $\pi(\hat M_k)=M$.
Then we get an induced surjective homomorphism
$\tilde\pi:\hat M_{k}/\hat M_{k-1}\longrightarrow M$.
So $M$ is a $\g'$-quotient of $\hat M_{k}/\hat M_{k-1}$,
which is a MOIS over $\g'$ by Theorem \ref{thm3.7}.
So is $M$, proving Theorem \ref{thm3.1}.

\section{Modules of intermediate series}
\def\theequation{4.\arabic{equation}}
\setcounter{equation}{0}

In this section we classify the irreducible modules of intermediate series for the algebra $\g$.
Through out this section we denote by $\overline{r}$ the residue of the integer $r$ by $p$.

Before we give the construction of MOIS for the algebra $\g$,
we recall such construction for the Virasoro algebra $Vir$, which appears in many references,
such as \cite{SZ}.
There are three kinds of MOIS over $Vir$,
denoted by $A(a),\ B(a),\ V(\al,\be)$ with parameters $a,\al,\be\in\C$,
which share a same basis
$\{w_i\mid i\in\Z\}$, and have $Vir$-actions as follows.\\
The action on $A(a)$:
$$K_xw_j=0,\ \  x_iw_j=\begin{cases}
                       (i+j)w_{i+j}, & \text{ if }j\neq 0;\\
                       i(i+a)w_i,  & \text{ if }j= 0.
                       \end{cases}$$
The action on $B(a)$:
$$K_xw_j=0,\ \  x_iw_j=\begin{cases}
                       jw_{i+j}, & \text{ if }i+j\neq 0;\\
                       -i(i+a)w_0,  & \text{ if }i+j= 0.
                       \end{cases}$$
The action on $V(\al,\be)$:
$$K_xw_j=0,\ \  x_iw_j=(\al+j+i\be)w_{i+j}.$$
The $Vir$-modules $A(a), B(a)$ are always reducible,
and $V(\al,\be)$ is reducible if and only if $\al\in\Z$ and $\be\in\{0,1\}$.
The $Vir$-module $V(\al,\be)$ with $\al\in \Z,\ \be\in\{0,1\}$ has a unique subquotient
denoted by $V'(\al,\be)$,
$A(a)$ has a unique subquotient $A'(a)$ isomorphic to $V'(0,1)$,
and $B(a)$ has a unique subquotient $B'(a)$ isomorphic to $V'(0,0)$.
Therefore in some sense $A(a)$ and $B(a)$ may be considered
as some "mutations" of $V(0,1)$ and $V(0,0)$.

Since the subalgebra $\gz$ of $\g$ is isomorphic to $Vir$,
and will play a key role in the rest of this section,
we turn the MOIS's over $Vir$ into the corresponding $\gz$-versions.
There are three kinds of MOIS over the algebra $\gz$,
denoted by $A_j(a),\ B_j(a),\ V_j(\al,\be)$ with parameters $a,\al,\be\in\C$,
which share a same basis
$\{\ve{j+pk}\mid k\in\Z\}$ (the index $j$, with $0\leq j\leq p-1$, is "redundant",
we put it here just for later narration convenience),
and respectively have $\gz$-action (let $m\in \pz$)
\begin{flalign*}
&\text{on }A_j(a):\ \ C_0\ve{j+pk}=0,\ \  L_m\ve{j+pk}=\begin{cases}
                       (m+pk)\ve{j+m+pk}, & \text{ if }k\neq 0;\\
                       m(m+a)v_{j+m},     & \text{ if }k= 0,
                       \end{cases}&\\
&\text{on }B_j(a):\ \ C_0\ve{j+pk}=0,\ \  L_m\ve{j+pk}=\begin{cases}
                       pkv_{j+m+pk},   & \text{ if }m+pk\neq 0;\\
                       -m(m+a)\ve{j},  & \text{ if }m+pk= 0,
                       \end{cases}&\\
&\text{and on }V_j(\al,\be):\ \ C_0\ve{j+pk}=0,\ \  L_m\ve{j+pk}=(\al+j+pk+m\be)v_{j+m+pk}.&
\end{flalign*}
We call these modules, or their subquotients {\em of type $A,B,V$} respectively.
The following lemma is obvious.

\begin{lem}\label{lem4.1}
(1) $V_{j+1}(\al,\be)\cong V_j(\al+1,\be)$.\\
(2) The $\gz$-module $V_j(\al,\be)$ is reducible
if and only if $\al\in -j+\pz$ and $\be\in\{0,1\}$.
Moreover, for $l\in\Z$, $V_j(-j-pl,1)$ has a unique subquotient(actually a submodule)
$$V'_j(-j-pl,1)=\spanc{\ve{-j+pk}\mid k\neq l},$$
and $V_j(-j-pl,0)$ has a unique quotient
$$V'_j(-j-pl,0)=V_j(-j-pl,0)/\C\ve{j+pl}.$$
(3) The $\gz$-module $A(a)$ is reducible and has a unique subquotient $A'_j(a)\cong V'_j(-j,1)$.\\
(4) The $\gz$-module $B(a)$ is reducible and has a unique subquotient $B'_j(a)\cong V'_j(-j,0)$.
\end{lem}

Now we give the construction of MOIS's over the algebra $\g$.
Let $\al,\be\in\C$ and
$F=(\f ij)$ be a $(p-1)\times p$ complex matrix,
with index $1\leq i\leq p-1,\ 0\leq j\leq p-1$,
satisfying the following three conditions
\begin{description}
\item[(I)] $0\in o(F)=\{j\mid \f ij\neq 0\text{ for some } i\}$;
\item[(II)] if $\f ij\neq 0$ then $\f s {\ov{i+j}}\neq 0$ for some $1\leq s\leq p-1$;
\item[(III)] $\f r{\ov{i+s}}\f si=\f s{\ov{i+r}}\f ri$
                for any $0\leq i\leq p-1$ and $1\leq r,s\leq p-1$;
\end{description}

For any $j\in o(F)$ denote $V_{(j)}=\spanc{\ve{j+pk}\mid k\in\Z}$.
Define the $\g$-module structure on $\V=\bigoplus_{j\in o(F)}V_{(j)}$ by
\begin{equation}\label{eq4.01}
\begin{aligned}
&L_m\ve{j+n}=(\al+j+n+m\be)\ve{j+n+m};\\
&L_s\ve{j+n}=\f{\ov s}j\ve{j+n+s};\\
&C_i\ve{j+n}=0,\ \ i=0,1,\cdots,p-1,
\end{aligned}
\end{equation}
where $m,n\in \pz,\ s\notin \pz$ and $j\in o(F)$.
Clearly $\V$ is a MOIS over $\g$ and $\V=V_0(\al,\be)$ if $o(F)=\{0\}$.
We call $V_{(j)}$ a {\em component} of $\V$,
call the component $V_{(i)}$ {\em directly links to} $V_{(j)}$ if $L_{j-i}V_{(i)}\neq 0$,
and call $V_{(i)}$ {\em links to} $V_{(j)}$
if they are directly linked through some other components.
The condition (II) makes sure that
each component links to some other component by some $L_s, s\notin \pz$.
Notice that each component $V_{(j)}$ of $\V$ is a MOIS over the algebra $\gz$.
This implies by Lemma \ref{lem4.1} that
if the order of $o(F)$ is more than 1, then the $\g$-module $\V$ must be irreducible.
Moreover we have

\begin{prop}\label{prop4.2}
(1) The $\g$-module $\V$ is reducible if and only if $o(F)=\{0\}$,
$\al\in \pz$ and $\be\in\{0,1\}$.\\
(2) The $\g$-modules $\V$ and $V(\al',\be',F')$ are isomorphic if and only if
$k=\al'-\al\in\Z, \be'=\be$ and $F'=\sgm^k(F\cdot D)$,
where $D=(D_{s,j})$ is a $(p-1)\times p$-matrix with all $D_{s,j}$ being nonzero,
$F\cdot D$ is defined to be the $(p-1)\times p$-matrix
with the $(s,j)$-th entry being $\f sj D_{s,j}$,
and $\sgm$ denotes the permutation $(0,1,2,\cdots,p-1)$
acting on a $(p-1)\times p$-matrix by shifting its columns.
\end{prop}
\begin{pf}
(1) is clear from the irreducibility of the $\gz$-module $V_0(\al,\be)$.\\
(2) Let $\phi:\V\longrightarrow V(\al',\be',F')$ be an isomorphism of $\g$-modules.
Then $\phi(V_{(0)})$ is a component of $V(\al',\be',F')$,
say $V'_{(j)}=\phi(V_{(0)})$ for some $j\in\{0,1,2,\cdots,p-1\}$,
that is, $V_{(0)}\cong V'_{(j)}$ as $\gz$-modules,
which implies $\be'=\be$ and $\al'+j-\al\in p\Z$.
Then $k=\al'-\al\in -j+\pz\subseteq\Z$.
Here, to denote notations of $V(\al',\be',F')$ we use the same symbols as $\V$
with an extra apostrophe.
Moreover we have $V_{(i)}\cong V'_{(i+j)}=V'_{(i-k)}$ as $\gz$-modules for any $i$.
These isomorphism maps are given by restrictions of $\phi$.
Without loss of generality we may assume
$$\phi(v_{i+n})=D_i v'_{i+n+k}, \text{ for }i\in o(F), n\in\pz
  \text{ and some }D_i\in\C\setminus\{0\}.$$
Now for any $i\in\Z, r\notin\pz$, since
$$\begin{aligned}
 \phi(L_rv_i)=D_{\ov {i+r}}\f{\ov r}{\ov i}v'_{i+r+k}\text{ and }
 L_r\phi(v_i)=D_{\ov i}F'_{{\ov r},{\ov{i+k}}}v'_{i+r+k},
\end{aligned}$$
we obtain that
$D_{\ov{i+r}}\f{\ov r}{\ov i}=D_{\ov i}F'_{{\ov r},{\ov{i+k}}}.$
Set $D=(D_{s,l})_{1\leq s\leq p-1,0\leq l\leq p-1}$ where
$D_{s,l}=\frac{D_{s+l}}{D_l}$. This proves $F'=\sgm^k(F\cdot D)$.

Conversely, it is easy to check that $v_i\mapsto D_{\ov i} v'_{i-k}$
defines a $\g$-module isomorphism from $\V$ to $V(\al',\be',F')$.
\pfend

The main result in this section is the following
\begin{thm}\label{thm4.3}
Let $V$ be an irreducible $\g$-module of intermediate series,
then $V$ must be the irreducible subquotient of $\V$
for some $\al,\be\in \C$ and $(p-1)\times p$ matrix $F$ satisfying conditions (I)-(III).
\end{thm}

From now on we fix $M$ an irreducible $\g$-module of intermediate series.
By the discussion about cuspidal modules in Section 3 we know that
$$C_iM=0\text{ for all }0\leq i\leq p-1.$$
Let $M=\bigoplus\limits_{i\in\Z}M_i$ be the weight space decomposition of $M$,
where $M_i$ is the weight space corresponding to
the weight $\al+i$ for a unanimous parameter $\al\in\C$, and $\dim M_i\leq 1$.
We note that some of these $M_i$'s may be $0$.
For $0\leq i\leq p-1$ set
$$M(i)=\bigoplus_{k\in\Z}M_{i+pk}.$$
We call $M(i)$, if not zero, a {\em component} of the $\g$-module $M$.
Clearly each component $M(i)$ forms a MOIS over $\gz$,
hence must be of type $A,B,V$.
Recall from Lemma \ref{lem4.1}
that the $\gz$-modules $A'_j(a)\cong V_j(-j,1)$ and $B'_j(a)\cong V'_j(-j,0)$.
To unify the $\gz$-actions on components $M(i)$ of $M$,
for $m\in \pz$ we rewrite the $L_m$-action on $M(i$)
$$L_m\ve{i+pk}=\begin{cases}
(\al_i+i+pk+m)\ve{i+pk+m}   &\text{ if $M(i)$ is of type }A,\ k\in\Z,\ k\neq 0;\\
(\al_i+i+pk)\ve{i+pk+m}     &\text{ if $M(i)$ is of type }B,\ k\in\Z,\ pk+m\neq 0,
\end{cases}$$
where $\al_i=-i$ actually.

Denote $o(M)=\{i\mid 0\leq i\leq p-1,\ M(i)\neq 0\}$,
and we have $M=\bigoplus_{i\in o(M)}M(i)$.
We may assume that $0$ lies in $o(M)$ by shifting the parameter $\al$ if necessary.
So if $M$ contains only one component, which forces $M=M(0)$,
then $L_sM=0$ for any $s\notin \pz$,
which means $M$ is an irreducible $\gz$-module,
hence isomorphic to $V_0(\al,\be)$, or its subquotient.
This proves Theorem \ref{thm4.3} for the case $o(M)=\{0\}$.

From now on we assume that $M$ contains more than one components.
Let $i\neq j\in o(M)$, we say that $M(i)$ {\em directly links to} $M(j)$ if $L_{j-i}M(i)\neq 0$,
and denote by $i\rightsquigarrow j$.
We say that $M(i)$ {\em links to }$M(j)$ if they are directly linked through some other components.
Since $M$ is irreducible, all components are linked in some way,
that is, for each $i\in o(M)$ there exist $j,k\in o(M)$
such that $j\rightsquigarrow i\rightsquigarrow k$.
This linkage will play an important part in the following proof of Theorem \ref{thm4.3}.

For $i\rightsquigarrow i+s\in o(M)$,
by Lemma \ref{lem4.1} we may write the $\g(s)$-action on $M(i)$ as follows.
$$L_{s+pl}\ve{i+pk}=\f {s+pl}{i+pk}\ve{i+s+pk+pl},\text{ for }k,l\in\Z,$$
where $\f {s+pl}{i+pk}\in\C$.
The fact that all components of $M$ are modules of intermediate series over $\gz$
and link to each other in some way
implies that all weight spaces $M_{j+pk}\neq 0$ if $j\in o(M), k\in\Z$,
and all $\f {s+pl}{i+pk}$ are nonzero.

Notice that the element $L_0$ has a same action on all three types $A,B,V$ of
MOIS's over $\gz$.
Fix $i\rightsquigarrow i+s\in o(M)$. Choose $k\neq 0,\pm1$ and we have
$$\begin{aligned}
s\f s{i+pk}\ve{i+s+pk}
&=sL_s\ve{i+pk}=[L_0,L_s]\ve{i+pk}=L_0L_s\ve{i+pk}-L_sL_0\ve{i+pk}\\
&=\big(\f s{i+pk}(\al_{i+s}+i+s+pk)-\f s{i+pk}(\al_i+i+pk)\big)\ve{i+s+pk},
\end{aligned}$$
where $\al_{i}=-i$ and $\al_{i+s}=-i-s$ if $M(i), M(i+s)$ are of type $A$ or $B$.
It follows that $\f s{i+pk}(\al_{i+s}-\al_i)=0$.
Since $\f s{i+pk}\neq 0$ by the choice of $k$, we have
$$\al_{i+s}=\al_i \text{ if }i\rightsquigarrow i+s.$$
This means that for all components $M(i)$ the parameters $\al_i$ coincide.
It follows that in the set $C(M)=\{M(i)\mid i\in o(V)\}$ of all components of $M$
there is at most one $M(i)$ that is of type $A$ or $B$.
So the set $C(M)$ must be one of the following three cases:\\
\indent (1) all $M(i)\in C(M)$ are of type $V$;\\
\indent (2) all $M(i)\in C(M)$ are of type $V$ except one $A$;\\
\indent (3) all $M(i)\in C(M)$ are of type $V$ except one $B$.\\
In the following we deal with these three cases separately.

\subsection{All components being of type $V$}

In this subsection we discuss the case where all components $M(i)$ are of type $V$.
Since for all components $M(i)$ the parameters $\al_i$ coincide,
we assume $M(i)=V_i(\al,\be_i)$ for some $\be_i\in\C$.
Therefore, fix $i,i+s\in o(M)$ such that $i\rightsquigarrow i+s$ and we have
$$\begin{aligned}
L_m\ve{i+pk}&=(\al+i+pk+\be_im)\ve{i+pk+m}\  \text{ for }m\in \pz, k\in\Z;\\
L_s\ve{i+pk}&=\f s{i+pk}\ve{i+pk+s}       \ \ \ \text{ for }s\notin \pz.
\end{aligned}$$

Let $m,n\in \pz$. Consider $[L_m,L_s]\ve{i+n}$ and we get
\begin{equation}\label{eq4.02}
s\f{m+s}{i+n}=\f s{i+n}(\al+i+n+s+\be_{i+s}m)-\f s{i+m+n}(\al+i+n+\be_im).
\end{equation}
Take $n=0$ in (\ref{eq4.02}) and we have
\begin{equation}\label{eq4.03}
\f s{i+m}(\al+i+\be_im)=\f s{i}(\al+i+s+\be_{i+s}m)-s\f{m+s}{i}.
\end{equation}
Consider for $k\in\Z$ that
$$\begin{aligned}
 s&(\al+i+\be_ikm)\f{m+s}{i+km}\ve{i+s+m+km}=(\al+i+\be_ikm)sL_{m+s}\ve{i+km}=sL_{m+s}L_{km}\ve i\\
  &=s[L_{m+s},L_{km}]\ve i+L_{km}L_{m+s}\ve i\\
  &=-s(m+s)L_{s+m+km}\ve i+s\f {m+s}i L_{km} \ve{i+s+m}\\
  &=s\big(-(m+s)\f{s+m+km}i+\f{m+s}i(\al+i+s+m+\be_{i+s}km)\big)\ve{i+s+m+km}.
\end{aligned}$$
On the other hand by (\ref{eq4.02}) (with $n$ replaced by $km$) and
(\ref{eq4.03})(with $m$ replaced by $km$) we have
$$\begin{aligned}
  &s(\al+i+\be_ikm)\f{m+s}{i+km}\\
 =&(\al+i+\be_ikm)\big(\f s{i+km}(\al+i+s+km+\be_{i+s}m)-\f s{i+m+km}(\al+i+km+\be_im)\big)\\
 =&\f si(\al+i+s+km+\be_{i+s}m)(\al+i+s+\be_{i+s}km)-s\f{km+s}i(\al+i+s+km+\be_{i+s}m)\\
  &\ \ -\f s{i+m+km}(\al+i+km+\be_{i}km)(\al+i+km+\be_{i}m).
\end{aligned}$$
Together we get
\begin{equation}\label{eq4.04}
\begin{aligned}
&s(m+s)\f{s+m+km}i+\f si(\al+i+s+km+\be_{i+s}m)(\al+i+s+\be_{i+s}km)\\
&-s\f{m+s}i(\al+i+s+m+\be_{i+s}km)-s\f{km+s}i(\al+i+s+km+\be_{i+s}m)\\
&=\f s{i+m+km}(\al+i+\be_{i}km)(\al+i+km+\be_{i}m).
\end{aligned}
\end{equation}
Rewrite (\ref{eq4.03}) as
$s\f{m+s}{i}=\f s{i}(\al+i+s+\be_{i+s}m)-\f s{i+m}(\al+i+\be_im).$
Substituting it into (\ref{eq4.04}) we get
\begin{equation}\label{eq4.05}
 \begin{aligned}
  &\f s{i+m}(\al+i+\bei m)(\al+i+s+m+\beis km)\\
  &+\f s{i+km}(\al+i+\bei km)(\al+i+s+km+\beis m)\\
  &-\f s{i+(k+1)m}(\al+i+\bei km)(\al+i+km+\bei m)\\
  &-\f s{i+(k+1)m}(m+s)\big(\al+i+\bei(k+1)m\big)\\
  =&\f si (\al+i+s+m+\beis km)(\al+i+s+\beis m)\\
  &-\f si(m+s)\big(\al+i+s+\beis(k+1)m\big) \text{ for any }k\in\Z, m\in\pz.
 \end{aligned}
\end{equation}
This recursive relation implies that $\f s{i+m}$ is a rational function in variable $m\in \pz$.
Write $\f s{i+m}=\frac{f(m)}{g(m)}$
where $f(m), g(m)$ are polynomials,
and we suppose $g(m)$ is not a constant.
Since the right hand side of (\ref{eq4.05}) is a polynomial,
to make (\ref{eq4.05}) valid the polynomial $g(km), k\in\Z$,
appearing in the denominator of the left hand side,
has to be cancelled out by the numerator obtained by reducing fraction to a common denominator.
But we can always replace $i$ by $i+n$ in (\ref{eq4.05}) with any $n\in\pz$.
This replacement would change the numerator in (\ref{eq4.05}),
which should be able to cancel out $g(km)$,
to a polynomial that is not proportional to the original one
and hence can not cancel out $g(km)$ anymore.
This deduces a contradiction.
So $g(m)$ is a constant and $\f s{i+m}$ is a polynomial.
Then by (\ref{eq4.03}) we see that $\f {m+s}i$ is a polynomial too.

Set $k=-1$ in (\ref{eq4.04}), and we get
\begin{equation}\label{eq4.06}
 \begin{aligned}
  &s\f{m+s}i(\al+i+s+m-\be_{i+s}m)+s\f{-m+s}i(\al+i+s-m+\be_{i+s}m)\\
  &=\f si\big((\be_i-\be_{i+s})(\be_i+\be_{i+s}-1)m^2+2s(\al+i+s)\big).
 \end{aligned}
\end{equation}
Let
\begin{equation}\label{eq4.07}
\f {m+s}i=\sum_{j=0}^ta_jm^j,
\end{equation}
where $t=\partial\f{m+s}i$ is the order of $\f{m+s}i$,
all $a_j\in\C,\ a_t\neq 0$ and $a_0=\f si\neq 0$.
Then the equation (\ref{eq4.06}) implies that
\begin{equation}\label{eq4.08}
\begin{aligned}
&(\al+i+s)a_{2k}+(1-\be_{i+s})a_{2k-1}=0\ \ \ \text{ for }k\geq 2,\\
&2s(\al+i+s)a_2+2s(1-\be_{i+s})a_1=a_0(\be_i-\be_{i+s})(\be_i+\be_{i+s}-1).
\end{aligned}
\end{equation}
Multiplying $(\al+i+\be_im+\be_ikm)$ to (\ref{eq4.04}) and using (\ref{eq4.03}) we get
\begin{equation}\label{eq4.09}
 \begin{aligned}
  &s\f{s+m+km}i(m+s)(\al+i+\be_i(k+1)m)\\
  &+s\f{s+m+km}i(\al+i+\be_ikm)(\al+i+km+\be_im)\\
  &-s\f{m+s}i(\al+i+\be_i(k+1)m)((\al+i+s+m+\be_{i+s}km))\\
  &-s\f{km+s}i(\al+i+\be_i(k+1)m)(\al+i+s+km+\be_{i+s}m)\\
  &+\f si(\al+i+\be_i(k+1)m)(\al+i+s+km+\be_{i+s}m)(\al+i+s+\be_{i+s}km)\\
  &-\f si(\al+i+\be_ikm)(\al+i+km+\be_{i}m)(\al+i+s+\be_{i+s}(k+1)m)=0.
 \end{aligned}
\end{equation}
Take $k=1$,
\begin{equation}\label{eq4.10}
\begin{aligned}
&2s(\al+i+2\be_im)(\al+i+s+m+\be_{i+s}m)\f{m+s}i\\
&-s(m+s)(\al+i+2\be_im)\f{2m+s}i-s(\al+i+\be_im)(\al+i+m+\be_im)\f{2m+s}i\\
&=\f si\big((\al+i+2\be_im)(\al+i+s+m+\be_{i+s}m)(\al+i+s+\be_{i+s}m)\\
&\ \ \ \ \ \ \ \ \ \ \ -(\al+i+\be_im)(\al+i+m+\be_{i}m)(\al+i+s+2\be_{i}m)\big).
\end{aligned}
\end{equation}

In the next we prove that $\f{m+s}i$ is actually independent of the variable $m$
in the following two key lemmas,
which deal with the case $\al+i+s=0$ in Lemma \ref{lem4.4}
and the case $\al+i+s\neq0$ in Lemma \ref{lem4.5} separately.
\begin{lem}\label{lem4.4}
Let $i\rightsquigarrow i+s$ and suppose that all $M(j)\in o(M)$ are of type $V$.
If $\al+i+s=0$, then $\be_{i+s}=\be_i$ and $\f {m+s}i=\f si$ for any $m\in \pz$.
\end{lem}
\begin{pf}
Note that $\al+i=-s\neq0$.
Substitute (\ref{eq4.07}) into (\ref{eq4.10}) and we get
\begin{equation}\label{eq4.11}
\begin{aligned}
 &2s(-s+2\be_im)(1+\be_{i+s})m(a_0+a_1m+a_2m^2+\cdots+a_tm^t)\\
-&s(-s+\be_im)(-s+m+\be_im)(a_0+2a_1m+4a_2m^2+\cdots+a_t2^tm^t)\\
-&s(m+s)(-s+2\be_im)(a_0+2a_1m+4a_2m^2+\cdots+a_t2^tm^t)\\
-&a_0(-s+2\be_im)(1+\be_{i+s})\be_{i+s}m^2+a_0(-s+\be_im)(-s+m+\be_{i}m)2\be_{i}m=0,
\end{aligned}
\end{equation}
where the coefficient of $m$ is $2a_0s^2(\be_{i+s}-\be_i)=0$, forcing
$$\be_{i+s}=\be_i.$$
Therefore equation (\ref{eq4.06}) becomes
\begin{equation}\label{eq4.12}
(\be_{i+s}-1)(\f{m+s}i-\f{-m+s}i)=0.
\end{equation}

{\em Case a: $\be_{i+s}=\be_{i}=1$}.
Suppose $t=\partial\f {m+s}i\geq 2$.
Then the leading term $m^{t+2}$ in (\ref{eq4.11}) has coefficient
$$-sa_t\be_i\big(\be_i2^t-4\be_{i+s}+3\times2^t-4\big)=-4sa_t(2^t-2)\neq 0,$$
which makes the left side of (\ref{eq4.11}) nonzero for $m$ large enough, a contradiction.
Hence $t\leq 1$ and equation (\ref{eq4.03}) becomes
$$\f s{i+m}(-s+m)=(a_0-sa_1)m-sa_0.$$
Noticing that $\f s{i+m}$ is a polynomial in $m$, we get
$\f s{i+m}=a_0=\f si$ and $a_1=0.$ So $\f{m+s}i=\f si$.

{\em Case b: $\be_{i+s}=\be_i\neq 1$}.
Then by (\ref{eq4.12}) we have $\f{m+s}i=\f{-m+s}i$ is even.
So we may write
$$\f{m+s}i=a_0+a_2m^2+\cdots+a_{t-2}m^{t-2}+a_tm^t,\ \text{with $t$ even.}$$
Then equation (\ref{eq4.10}) becomes
\begin{equation}\label{eq4.13}
\begin{aligned}
 &2s(-s+2\be_im)(1+\be_{i+s})m(a_0+a_2m^2+\cdots+a_{t-2}m^{t-2}+a_tm^t)\\
-&s(-s+\be_im)(-s+m+\be_im)(a_0+4a_2m^2+\cdots+a_{t-2}2^{t-2}m^{t-2}+a_t2^tm^t)\\
-&s(m+s)(-s+2\be_im)(a_0+4a_2m^2+\cdots+a_{t-2}2^{t-2}m^{t-2}+a_t2^tm^t)\\
-&a_0(-s+2\be_im)(1+\be_{i+s})\be_{i+s}m^2+a_0(-s+\be_im)(-s+m+\be_{i}m)2\be_{i}m=0.
\end{aligned}
\end{equation}
If $t>2$ then the term $m^{t+1}$ in (\ref{eq4.13}) has coefficient $2a_ts^2(2^t-\be_{i+s}-1)=0$,
forcing $\be_{i+s}=\be_i=2^t-1$.
Then the leading term $m^{t+2}$ in (\ref{eq4.13}) has coefficient
$$-sa_t\be_i(\be_i2^t-4\be_{i+s}+3\times2^t-4)=-sa_t2^t(2^t-1)(2^t-2)\neq 0.$$
This provides a contradiction, so $t\leq2$.

If $t=2$ then the coefficient of $m^4$ in (\ref{eq4.13}) is $-8sa_2\be_i=0$,
forcing $\be_i=0$.
Hence the left hand side of equation (\ref{eq4.13}) turns to
$$-2s^2m(a_0+a_2m^2)+s^2(-s+m)(a_0+4a_2m^2)+s^2(m+s)(a_0+4a_2m^2)
  =6a_2s^2m^3\neq 0.$$
This is a contradiction. So $t=0$ and $\f {m+s}i=\f si$.
\pfend

Next we deal with the case $\al+i+s\neq0$.
\begin{lem}\label{lem4.5}
Let $i\rightsquigarrow i+s$ and suppose all $M(j)\in o(M)$ are of type $V$.
If $\al+i+s\neq0$, then $\beis=\bei$ and $\f{m+s}i=\f si$ for any $m\in \pz$.
\end{lem}
\begin{pf}
Suppose $\al+i+s\neq0$. Substitute (\ref{eq4.07}) into (\ref{eq4.10}) and we get
\begin{equation}\label{eq4.14}
\begin{aligned}
 &2s(\al+i+2\be_im)(\al+i+s+m+\be_{i+s}m)(a_0+a_1m+a_2m^2+\cdots+a_tm^t)\\
+&s(\al+i+\be_im)(\al+i+m+\be_im)(a_0+2a_1m+4a_2m^2+\cdots+a_t2^tm^t)\\
-&s(m+s)(\al+i+2\be_im)(a_0+2a_1m+4a_2m^2+\cdots+a_t2^tm^t)\\
-&a_0(\al+i+2\be_im)(\al+i+s+m+\be_{i+s}m)(\al+i+s+\be_{i+s}m)\\
+&a_0(\al+i+\be_im)(\al+i+m+\be_{i}m)(\al+i+s+2\be_{i}m)=0,
\end{aligned}
\end{equation}
in which the coefficient of the term $m^{t+2}$ is
\begin{equation}\label{eq4.15}
sa_t\bei(2^t\bei+4\be_{i+s}-2^t+4)=0\text{ if }t\geq 2,
\end{equation}
the coefficient of $m$ is
\begin{equation}\label{eq4.16}
4sa_1(\al+i)^2+2a_0(\al+i)(2s\be_i+(\al+i)\be_i-(\al+i)\be_{i+s}+s)=0,
\end{equation}
and the coefficient of $m^2$ is
\begin{equation}\label{eq4.17}
\begin{aligned}
&a_0(2s\be_i(1+\be_i)+(\al+i)(5\be_i^2+\bei-\be_{i+s}^2-\beis-4\bei\beis))\\
&+2sa_1(\al+i)(\be_{i+s}+4\be_i+1)+2sa_2(\al+i)(3(\al+i)-s)=0.
\end{aligned}
\end{equation}
On the other hand, replace $s$ by $m+s$ in equation (\ref{eq4.06}) and we get
\begin{equation}\label{eq4.18}
\begin{aligned}
&(m+s)(\al+i+s+2m-\beis m)\f{2m+s}i+(m+s)(\al+i+s+\beis m)\f si\\
&\ \ \ \ -\f{m+s}i\big((\bei-\beis)(\bei+\beis-1)m^2+2(m+s)(\al+i+s+m)\big)=0,
\end{aligned}
\end{equation}
that is

$$\begin{aligned}
&(m+s)(\al+i+s+2m-\beis m)(a_0+2a_1m+4a_2m^2+\cdots+a_t2^tm^t)\\
&-(\bei-\beis)(\bei+\beis-1)m^2(a_0+a_1m+a_2m^2+\cdots+a_tm^t)\\
&-2(m+s)(\al+i+s+m)(a_0+a_1m+a_2m^2+\cdots+a_tm^t)\\
&+a_0(m+s)(\al+i+s+\beis m)=0,
\end{aligned}$$
where the coefficient of $a_t m^{t+2}$ is
\begin{equation}\label{eq4.19}
\beis^2-\beis(2^t+1)+2^{t+1}-2-\bei^2+\bei=0\text{ if }t\geq2,
\end{equation}
and the coefficient of $m^2$ is
\begin{equation}\label{eq4.20}
2sa_1(1-\beis)-a_0(\bei-\beis)(\bei+\beis-1)+2sa_2(\al+i+s)=0.
\end{equation}
To prove our lemma we make\\
\indent{\em Claim 1: If $\al+i+s\neq0$ and $\al+i\neq0$ then $t=\partial\f{m+s}i<2$.}\\
\indent{\em Claim 2: If $\al+i+s\neq0$ and $\al+i\neq0$ then $\f{m+s}i=\f si$.
            Moreover, we have $\beis=\bei$ except for one case where $\bei=0$ and $\beis=1$.}\\
\indent{\em Claim 3: If $\al+i+s\neq0$ and $\al+i=0$ then $\f{m+s}i=\f si$ and $\beis=\bei$.}\\
\indent{\em Claim 4: The case where $\al+i+s\neq0,\ \al+i\neq0,\ \bei=0$ and $\beis=1$
            is self-contradictory.}\\
Clearly the lemma follows from these four claims.

{\em Proof of Claim 1}:
Suppose $t\geq 2$. Then by (\ref{eq4.15}) we get
$$\bei=0\text{ or }\beis=2^{t-2}(1-\bei)-1.$$
We prove Claim 1 in the following two cases.

{\em Case 1.1: $\al+i\neq 0$ and $\bei=0$}.
The equation (\ref{eq4.19}) implies that
$$\beis=2\text{ or }\beis=2^t-1.$$
Replacing $s$ by $s+km\ (k\in\Z)$ in (\ref{eq4.06}) we get
\begin{equation}\label{eq4.21}
\begin{aligned}
&\f{s+m+km}i(s+km)(\al+i+s+m+km-\beis m)\\
&+\f{s+km-m}i(s+km)(\al+i+s+km-m+\beis m)\\
&-\f{s+km}i\big((\bei-\beis)(\bei+\beis-1)m^2+2(s+km)(\al+i+s+km)\big)=0,
\end{aligned}
\end{equation}
in which the coefficient of $a_tm^{t+2}$ is
\begin{equation}\label{eq4.22}
k(k+1-\beis)(k+1)^t+(k-1)^tk(k-1+\beis)-k^t(\bei-\beis)(\bei+\beis-1)-2k^{t+2}=0.
\end{equation}
If $\beis=2$ then the left side of (\ref{eq4.22}) is
$$k^t(t-1)(t-2)+\cdots,$$
and if $\beis=2^t-1$ then the left side of (\ref{eq4.22}) is
$$k^t(2^t-t-1)(2^t-t-2)+\cdots,$$
respectively,
which are both nonzero for $t>2$, contradicting to (\ref{eq4.22}) with sufficiently large $k$.
So $t=2$ and then $\beis=2$ or $\beis=3$.

If $t=2$ and $\beis=2$, then the left side of (\ref{eq4.21}) turns to
$$2(\al+i)a_2km^3+2m^2(s^2a_2+sa_2(\al+i)-sa_1+a_0),$$
which is nonzero if $k$ and $m$ are chosen large enough, a contradiction.

If $t=2$ and $\beis=3$, then the equations (\ref{eq4.16}), (\ref{eq4.17}) and (\ref{eq4.20})
give a system of homogeneous linear equations on variables $a_0,a_1,a_2$
$$\left(\begin{array}{ccc}
   -4s        &      6       & 2s(\al+i+s)\\
2s(\al+i)    &  s-3(\al+i)  & 0\\
8s(\al+i)    &  -12(\al+i)  & 2s(\al+i)(3(\al+i)-s)\\
\end{array}\right)
\left(\begin{array}{c}
a_1\\a_0\\a_2
\end{array}\right)=0,
$$
which has a nonzero solution with $a_0\neq0$ and $a_2\neq0$.
It follows that
\begin{equation}\label{eq4.23}
s=-5(\al+i)\text{ and }a_0=4a_2(\al+i)^2.
\end{equation}
Moreover, substitute (\ref{eq4.07}) into (\ref{eq4.21}) and we get
$$2km^3(a_2(\al+i-2s)+a_1)+2m^2(sa_2(\al+i+s)-2sa_1+3a_0)=0.$$
It follows by (\ref{eq4.23}) that
$$a_0=30a_2(\al+i)^2,$$
which implies $a_2=0$,
contradicting to $t=2$. So Claim 1 stands in Case 1.1.

{\em Case 1.2: $\al+i\neq0,\ \bei\neq0$ and $\beis=2^{t-2}(1-\bei)-1$}.
Consider (\ref{eq4.22}) as a polynomial in variable $k$.
Since it stands for any $k\in\Z$, the coefficient of $k^t$ is
$$\beis^2-\beis(2t+1)+t^2+t-\bei^2+\bei=0.$$
Subtracting equation (\ref{eq4.19}) gives
$$\beis (2^t-2t)=2^{t+1}-t^2-t-2.$$
Suppose $t>2$.
Since the coefficient of $k^{t-2}$ in (\ref{eq4.22}) is
$\frac{1}{6}\binom t2 (t-4\beis+1)=0$,
it follows
$$\beis=\frac{t+1}{4}.$$
So we have $(7-t)2^t=2t^2+2t+8$. Hence
$$t=3\text{ or }t=4.$$
If $t=3$, then $\beis=1$ and $\bei=0$, which is a contradiction.
If $t=4$, then $\beis=\frac{5}{4}$ and $\bei=\frac{7}{16}$.
So (\ref{eq4.22}) has leading term $\frac{2577}{256}k^4\neq 0$,
which is a contradiction if $k$ is large enough.

Suppose $t=2$ and we have $\beis=-\bei$.
Then (\ref{eq4.18}) becomes
$$\begin{aligned}
&6a_2k^2m^4(1+\bei)+2km^3[a_2(\al+i+2s)+2sa_2(1+\bei)+a_1(1+2\bei)]\\
&\ \ \ \ +2m^2[sa_2(\al+i+s)+sa_1(1+\bei)+a_0\bei]=0,
\end{aligned}$$
from which we have
$$\bei=-1,\beis=1, a_1=(\al+i+2s)a_2,\ a_0=s(\al+i+s)a_2.$$
Then equation (\ref{eq4.14}) turns to
$$\begin{aligned}
&2s(\al+i-2m)(\al+i+s+2m)(a_0+a_1m+a_2m^2)\\
&+s(\al+i-m)(\al+i)(a_0+2a_1m+4a_2m^2)\\
&-s(m+s)(\al+s-2m)(a_0+2a_1m+4a_2m^2)\\
&-a_0(\al+i-2m)(\al+i+s+m)(\al+i+s+2m)\\
&+a_0(\al+i-m)(\al+i)(\al+i+s-2m)=0,
\end{aligned}$$
where the coefficient of $m^3$ is $-8sa_2(\al+i)\neq0$,
which is a contradiction.
So Claim 1 stands in Case 1.2.

{\em Proof of Claim 2}:
Suppose otherwise, then by Claim 1 we have $t=1$.
Write $\f{m+s}i=a_0+a_1m$ where $a_1\neq0$ and $a_0=\f si\neq0$.
Then equation (\ref{eq4.14}) turns to
\begin{equation}\label{eq4.24}
\begin{aligned}
 &2s(\al+i+2\be_im)(\al+i+s+m+\be_{i+s}m)(a_0+a_1m)\\
+&s(\al+i+\be_im)(\al+i+m+\be_im)(a_0+2a_1m)\\
-&s(m+s)(\al+i+2\be_im)(a_0+2a_1m)\\
-&a_0(\al+i+2\be_im)(\al+i+s+m+\be_{i+s}m)(\al+i+s+\be_{i+s}m)\\
+&a_0(\al+i+\be_im)(\al+i+m+\be_{i}m)(\al+i+s+2\be_{i}m)=0,
\end{aligned}
\end{equation}
in which the coefficient of $m$ is
\begin{equation}\label{eq4.25}
4sa_1(\al+i)^2+2a_0(\al+i)((\al+i)(\bei-\beis)+2s\bei+s)=0,
\end{equation}
the coefficient of $m^2$ is
\begin{equation}\label{eq4.26}
2sa_1(\al+i)(\beis+4\bei+1)+a_0[2s\bei(\bei+1)+(\al+i)(5\bei^2+\bei-\beis^2-\beis-4\bei\beis)]=0,
\end{equation}
and the coefficient of $m^3$ is
\begin{equation}\label{eq4.27}
2sa_1\bei(2\beis+\bei+1)+2a_0\bei(\bei-\beis)(\bei+\beis+1)=0.
\end{equation}
Moreover, equation (\ref{eq4.18}) becomes
$$\begin{aligned}
&(m+s)(\al+i+s+2m-\beis m)(a_0+2a_1m)-(\bei-\beis)(\bei+\beis-1)m^2(a_0+a_1m)\\
&-2(m+s)(\al+i+s+m)(a_0+a_1m)+a_0(m+s)(\al+i+s+\beis m)=0,
\end{aligned}$$
in which the coefficient of $m^2$ is
\begin{equation}\label{eq4.28}
2sa_1(1-\beis)-a_0(\bei-\beis)(\bei+\beis-1)=0,
\end{equation}
and the coefficient of $m^3$ is
\begin{equation}\label{eq4.29}
a_0(\beis^2-3\beis-\bei^2+\bei+2)=0.
\end{equation}
We prove $\f{m+s}i=\f si$ in the following three cases.

{\em Case 2.1: $\al+i+s\neq 0,\ \al+i\neq 0$ and $\bei=0$.}
In this case we have $\beis=1$ or $\beis=2$ by (\ref{eq4.29}).
If $\beis=1$ then equations (\ref{eq4.25}) and (\ref{eq4.26})
give a system of homogeneous linear equations with variables $a_0, a_1$
$$\left(\begin{array}{cc}
2s(\al+i)  &  s-\al-i\\
4s(\al+i)  &  -2(\al+i)\\
\end{array}\right)
\left(\begin{array}{c}
a_1\\a_0
\end{array}\right)=0.$$
Since the determinant of its coefficient matrix is $-4s^2(\al+i)\neq0$,
it forces $a_0=a_1=0$,
which is a contradiction.

If $\beis=2$ then equations (\ref{eq4.25}) and (\ref{eq4.26})
give a system of homogeneous linear equations with variables $a_0, a_1$
$$\left(\begin{array}{cc}
2s(\al+i)  &  s-2\al-2i\\
6s(\al+i)  &  -6(\al+i)\\
\end{array}\right)
\left(\begin{array}{c}
a_1\\a_0
\end{array}\right)=0,$$
the determinant of whose coefficient matrix is $-6s^2(\al+i)\neq 0$,
also a contradiction.
So in Case 2.1 we have $\f{m+s}i=\f si$.

{\em Case 2.2: $\al+i+s\neq 0,\ \al+i\neq 0,\ \bei\neq 0$ and $\beis=\bei$.}
In this case equation (\ref{eq4.28}) becomes $2sa_1(1-\beis)=0$,
forcing $\bei=\beis=1$. Then equation (\ref{eq4.27}) becomes $8sa_1=0$,
which contradicts to $a_1\neq0$.
So $\f{m+s}i=\f si$ in Case 2.2.

{\em Case 2.3: $\al+i+s\neq 0,\ \al+i\neq 0,\ \bei\neq 0$ and $\beis\neq\bei$.}
Since $a_0\neq 0, a_1\neq 0$, the determinant of the coefficient matrix of the system of
homogeneous linear equations provided by (\ref{eq4.27}) and (\ref{eq4.28}) must be 0,
that is,
$$2\bei(\bei+1)(\beis+\bei+1)=0.$$
So $\bei=-1$ or $\beis+\bei+1=0$.

Suppose $\bei=-1$.
The fact that the determinant of the coefficient matrix of the system of
homogeneous linear equations given by (\ref{eq4.26}) and (\ref{eq4.27}) is 0 gives that
$$(\al+i)\beis(\beis+1)=\beis(\al+i+s+(\al+i)\beis),$$
forcing $\beis=0$. Then we have from (\ref{eq4.25}) and (\ref{eq4.28}) that
$$\left(\begin{array}{cc}
2s(\al+i)  &  -s-\al-i\\
    s      &    -1\\
\end{array}\right)
\left(\begin{array}{c}
a_1\\a_0
\end{array}\right)=0.$$
So $a_0=sa_1$ and $2s(\al+i)=s(\al+i+s)$, hence $s=\al+i$.
Moreover, from equation (\ref{eq4.26}) we have
$$0=-6s(\al+i)a_1+(4\al+4i-s)a_0=-6s^2a_1+3sa_0,$$
it follows that $a_0=2sa_1$, contradicting to $a_0=sa_1$.

Suppose that $\bei\neq -1$ and $\beis+\bei+1=0$.
Then equation (\ref{eq4.27}) becomes $sa_1\beis=0$, forcing $\beis=0$.
Therefore $\bei=-1-\beis=-1$, contradiction.
This validates $\f{m+s}i=\f si$ in Case 2.3.

Now equation (\ref{eq4.10}) turns to
$$2m^3\bei\beis(\beis-\bei)+m^2(\al+i)(\beis-\bei)(\beis+\bei-1)=0.$$
Suppose $\beis\neq\bei$ and $\bei\neq 0$,
then the equation above implies that $\beis=0$ and $\bei=1-\beis=1$.
Since $\f{m+s}i=\f si$, equation (\ref{eq4.03}) turns to
$$\f s{i+m}(\al+i+m)=\f si(\al+i),$$
which is a contradiction since
the left hand side is a polynomial of order at least 1,
but the right hand side is a constant.
So $\beis=\bei$, or $\bei=0$.
If $\beis\neq\bei$, then $\bei=0$ and the coefficient of $m^2$ being 0 forces
$\beis=1-\bei=1$, as in the exceptional case stated in Claim 2.

{\em Proof of Claim 3}:
Notice that $\al+i=0$ and then equation (\ref{eq4.17}) becomes
$$2sa_0\be_i(1+\be_i)=0,$$
which forces $\bei=0$ or $\bei=-1$.
We prove Claim 3 in two cases.

{\em Case 3.1: $\al+i+s\neq 0,\ \al+i=0$ and $\bei=0$}.
Let $n=0$ in (\ref{eq4.02}) and we have
$$\f{m+s}i=\frac{a_0}{s}(\beis m+s);\ \ t\leq1;\ \ a_1=\frac{a_0\beis}{s}.$$
Substitute into (\ref{eq4.18}) and we get
$$a_0m^3\beis(\beis-1)(\beis-2)-sa_0m^2\beis(\beis-1)=0.$$
So $\beis=0$ or $\beis=1$.
If $\beis=0=\bei$ then $\f{m+s}i=\f si$ as claimed.

Suppose $\beis=1$ and let $i+s\rightsquigarrow i+s+q\in o(M)$.
Notice that $\al+i+s\neq0$ and $\al+i=0$.
If $\al+i+s+q=0$, i.e. $s+q=0$, then $M(i+s+q)=M(i)$ and $\be_{i+s+q}=\bei=0$,
but by Lemma \ref{lem4.4} we have $\be_{i+s+q}=\beis=1$.
This contradiction forces $\al+i+s+q\neq0$.
Hence by Claim 2 we get $\be_{i+s+q}=\beis=1$.
Let $i+s+q\rightsquigarrow i+s+q+r$.
A similar proof to that of the $i+s\rightsquigarrow i+s+q$ case shows $\be_{i+s+q+r}=1$.
This proves that the linkage from $M(i+s)$ always goes to
components $M(j)$ with parameter $\be_j=1$, hence never back to $M(i)$,
which implies that $M$ is reducible, a contradiction.
This proves Claim 3 in Case 3.1.

{\em Case 3.2: $\al+i+s\neq 0,\ \al+i=0$ and $\bei=-1$}.
In this case (\ref{eq4.15}) turns to
$$-4sa_t(\beis-2^{t-1}+1)=0,$$
forcing $\beis=2^{t-1}-1$.
Then (\ref{eq4.19}) becomes
$$-(2^{t-1}-1)(2^{t-1}-2)=0,$$
which is absurd if $t>2$.
Suppose $t=2$, and then $\beis=1$.
Therefore (\ref{eq4.14}) turns to
$$-4sm(s+2m)(a_0+a_1m+a_2m^2)+2sm(m+s)(a_0+2a_1m+4a_2m^2)+2a_0m(s+m)(s+2m)=0,$$
where the coefficient of $m^3$ is $4a_0\neq 0$.
This is a contradiction.
Hence $t\leq 1$ and we write $\f{m+s}i=a_0+a_1m$.
The equation (\ref{eq4.10}) becomes
$$2m^3\big(a_0\beis(\beis+1)-2sa_1\beis\big)=0,$$
and (\ref{eq4.20}) becomes
$$a_0(\beis+1)(\beis-2)+2sa_0(1-\beis)=0.$$
Hence we get a system of homogeneous linear equations in $a_0,a_1$
\begin{equation}\label{eq4.30}
\left(\begin{array}{cc}
2s(1-\beis) & (\beis+1)(\beis-2)\\
-2s\beis    & \beis(\beis+1)\\
\end{array}\right)
\left(\begin{array}{c}
a_1\\a_0\\
\end{array}\right)=0,
\end{equation}
which has a solution with $a_0\neq 0$.
So the determinant of the coefficient matrix is $-2s\beis(\beis+1)=0$,
forcing $\beis=-1$ or $\beis=0$.

Suppose $\beis=0$.
Let $i+s\rightsquigarrow i+s+q\in o(M)$.
Since $\al+i+s\neq0$,
we see $\be_{i+s+q}=0$ or 1 by Lemma \ref{lem4.4} and Claim 2.
So $M(i+s+q)\neq M(i)$ and $\al+i+s+q\neq 0$ since $\bei=-1$.
Therefore if we let $i+s+q\rightsquigarrow i+s+q+r$,
then we still have $\be_{i+s+q+r}=0$ or 1 by Lemma \ref{lem4.4} and Claim 2.
Hence $M(i+s+q+r)\neq M(i)$ and $\al+i+s+q+r\neq 0$.
This means the linkage starting from $M(i)$
$$i\rightsquigarrow i+s\rightsquigarrow i+s+q\rightsquigarrow i+s+q+r\rightsquigarrow\cdots$$
always goes to component $M(j)$
with parameter $\be_j=0$ or 1, hence never back to $M(i)$,
which contradicts to the irreducibility of $M$.
This contradiction makes $\beis=-1=\bei$,
then the solution of the system of equations (\ref{eq4.30})
must be such that $a_1=0$.
So $\f{m+s}i=\f si$ as claimed.

Now equation (\ref{eq4.03}) turns to
$$\bei\f s{i+m}=\beis\f si.$$
This implies that $\bei=0$ is equivalent to $\beis=0$.
Let $\bei\neq0$, i.e. $\beis\neq0$.
The equation (\ref{eq4.04}) becomes
$$km^2\beis(\beis-\bei)=0,$$
which implies $\beis=\bei$, as in Claim 3.

{\em Proof of Claim 4:}
Let $i+s\rightsquigarrow i+s+q$.
We have $\be_{i+s+q}=\beis=1$ by Lemma \ref{lem4.4} if $\al+i+s+q=0$,
and by Claim 2 if $\al+i+s+q\neq0$.
Furthermore, if we let $i+s+q\rightsquigarrow i+s+q+r$,
then we get $\be_{i+s+q+r}=\be_{i+s+q}=1$ by Lemma \ref{lem4.4}, Claim 2 or Claim 3
depending on whether $\al+i+s+q$ and $\al+i+s+q+r$ are 0 or not, separately.
Similarly we have $\be_{i+s+q+r+l}=\be_{i+s+q+r}=1$ if $i+s+q+r\rightsquigarrow i+s+q+r+l$.
This linkage goes on through components $M(j)$ with parameters all being $\be_j=1$,
hence it never goes back to $M(i)$, which is a contradiction.
This proves Claim 4.
\pfend

Now since the parameter $\be_j$ for components $M(j)$ in $M$ coincide,
we write $\be=\be_j$ for any $j\in o(M)$.
Furthermore, we have
\begin{lem}\label{lem4.6}
Let $i\rightsquigarrow i+s$. Then we have $\f{m+s}{n+i}=\f si$ for any $m,n\in \pz$.
\end{lem}
\begin{pf}
First we claim that $\f{s}{n+i}=\f si$ for any $n\in \pz$.
Since $\beis=\bei=\be$ and $\f{m+s}i=\f si$ for any $m\in \pz$,
we have by (\ref{eq4.03})
$$(\f s{i+n}-\f si)(\al+i+\be n)=0,\text{ for any }n\in \pz,$$
which implies the claim unless $\al+i+\be n=0$.
If $\al+i\neq 0$ and there are nonzero $n\in \pz$ such that $\al+i+\be n=0$,
then we have $\be\neq 0$.
Since $\al+i+n+\be n=n\neq0$, we have
$$\f s{i+n}=\f s{i+n+n}=\f s{i+2n}.$$
On the other hand, $\al+i+2\be n\neq0$ implies that $\f s{i+2n}=\f si$.
So $\f s{i+n}=\f si$.

If $\al+i=0$ and there are nonzero $n\in \pz$ such that $\al+i+\be n=0$,
then we have $\be=0$.
Since $\al+i+n+\be(-n)=n\neq 0$, we have
$\f s{i+n}=\f s{i+n+(-n)}=\f si$ as claimed.
Now the lemma follows from (\ref{eq4.02}).
\pfend

Now by Lemma \ref{lem4.4}, Lemma \ref{lem4.5} and Lemma \ref{lem4.6},
we see that the $\g$-action on $M$ satisfies equation (\ref{eq4.01}).
Set $F=(\f ij)_{1\leq i\leq p-1,\ 0\leq j\leq p-1}$.
Clearly the conditions (I) and (II) follow from the irreducibility of $M$,
and condition (III) follows directly from the Lie bracket
$[L_r,L_s]=\dt_{r+s,0}rC_{\overline r}$ and $C_{\overline r}M=0$.
This proves Theorem \ref{thm4.3} if all components of $M$ are of type $V$.
Then to complete the proof of Theorem \ref{thm4.3} we only need to show that
the module $M$ in the cases (2) and (3) does not exist.

\subsection{All components being of type $V$ except one $A$}

In this subsection we deal with the case where all components of $M$ are
of type $V$ except one $A$.
We may assume $i\in o(M)$ and $M(i)=A_i(a)$.
By our unification of $\gz$-actions on components of $M$,
we see that in this case the unanimous parameter $\al$ equals to $-i$.

Let $i-r\in o(M)$ and $i-r\rightsquigarrow i$.
We may assume $M(i-r)=V_{i-r}(-i,\beir)$ for some $\beir\in\C$.
Consider for any $m\in \pz$ and $0\neq n\in \pz$
$$\begin{aligned}
&(n+r)\f {r+m+n}{i-r}\ve{i+m+n}=(n+r)L_{r+m+n}\ve{i-r}=[L_m,L_{n+r}]\ve{i-r}\\
&=(L_mL_{n+r}-L_{n+r}L_m)\ve{i-r}=\big(\f{n+r}{i-r}(m+n)-\f{n+r}{m+i-r}(-r+\beir m)\big)\ve{i+m+n}.
\end{aligned}$$
So
$$(n+r)\f {r+m+n}{i-r}=\f{n+r}{i-r}(m+n)-\f{n+r}{m+i-r}(-r+\beir m), \text{ for }n\neq0.$$
Replace $m$ by $km$ and $n$ by $m$ and we get
\begin{equation}\label{eq4.31}
\begin{aligned}
\f{m+r}{km+i-r}(-r+\beir km)=&\f{m+r}{i-r}(k+1)m-(m+r)\f {r+(k+1)m}{i-r},\\
                             &\text{ for }m\neq0 \text{ and any }k\in\Z.
\end{aligned}
\end{equation}
For any $m\in \pz$ consider
$$\begin{aligned}
 &r\f{m+r}{i-r}\ve{i+m}=rL_{m+r}\ve{i-r}=[L_m,L_r]\ve{i-r}\\
=&\big(\f r{i-r}m(m+a)-\f r{m+i-r}(-r+\beir m)\big)\ve{i+m},
\end{aligned}$$
which gives,
$$r\f{m+r}{i-r}=\f r{i-r}m(m+a)-\f r{m+i-r}(-r+\beir m),$$
or
\begin{equation}\label{eq4.32}
\f r{m+i-r}(-r+\beir m)=\f r{i-r}m(m+a)-r\f{m+r}{i-r}.
\end{equation}
Moreover, for $0\neq k\in\Z$ consider
$$\begin{aligned}
 &r\f{m+r}{i-r+km}\ve{i+m+km}=[L_m,L_r]\ve{i-r+km}\\
=&\big(\f r{i-r+km}(k+1)m-\f r{i-r+m+km}(-r+km+\beir m)\big)\ve{i+m+km}.
\end{aligned}$$
Hence
\begin{equation}\label{eq4.33}
r\f{m+r}{i-r+km}=\f r{i-r+km}(k+1)m-\f r{i-r+m+km}(-r+km+\beir m), \text{ for }k\neq0.
\end{equation}
Then by (\ref{eq4.31}) and (\ref{eq4.33}) we have for any $k\neq0$ and $m\ (\neq 0)\in \pz$
\begin{equation}\label{eq4.34}
\begin{aligned}
&r\f{m+r}{i-r}(k+1)m-r(m+r)\f {r+(k+1)m}{i-r}=r\f{m+r}{km+i-r}(-r+\beir km)\\
&=(-r+\beir km)\big(\f r{i-r+km}(k+1)m-\f r{i-r+m+km}(-r+km+\beir m)\big)
\text{ for }k\neq 0.
\end{aligned}
\end{equation}
Take $k=-1$ in (\ref{eq4.34}) and we get
$$\f r{i-r}\beir(1-\beir) m^2=0, \text{ for any }m\neq0,$$
which forces $\beir=0$ or $\beir=1$.

Multiplying (\ref{eq4.34}) by $(-r+\beir (k+1)m)$ and applying (\ref{eq4.32}) we get
$$\begin{aligned}
&r\f{m+r}{i-r}(k+1)m(-r+\beir (k+1)m)-r(m+r)\f{km+m+r}{i-r}(-r+\beir (k+1)m)\\
=&\f r{i-r}k(k+1)m^2(km+a)(-r+\beir (k+1)m)\\
&-(-r+km+\beir m)(-r+\beir km)\big(\f r{i-r}(k+1)m((k+1)m+a)-r\f{(k+1)m+r}{i-r}\big)\\
&-r(k+1)m\f{km+r}{i-r}(-r+\beir (k+1)m)\text{ for }k\neq 0,\\
\end{aligned}$$
from which we see that $\f {m+r}{i-r}$ is a polynomial in $m$
for a similar reason to that for $\f s{i+m}$ in Subsection 4.1.
Write
$$\f {m+r}{i-r}=\sum_{j=0}^ta_jm^j, $$
where $t$ denotes the order of this polynomial and $a_j\in\C, a_t\neq0, a_0=\f {r}{i-r}\neq0$.
So the equation above turns to
\begin{equation}\label{eq4.35}
\begin{aligned}
0=&r(k+1)m(-r+\beir (k+1)m)(a_0+a_1km+\cdots+a_tk^tm^t)\\
+&r(k+1)m(-r+\beir (k+1)m)(a_0+a_1m+\cdots+a_tm^t)\\
-&r(k+1)m(-r+\beir (k+1)m)(a_0+a_1(k+1)m+\cdots+a_t(k+1)^tm^t)\\
-&r(-r+km+\beir m)(-r+\beir km)(a_0+a_1(k+1)m+\cdots+a_t(k+1)^tm^t)\\
+&a_0(k+1)m((k+1)m+a)(-r+km+\beir m)(-r+\beir km)\\
-&a_0k(k+1)m^2(km+a)(-r+\beir (k+1)m)\text{ for } k\neq 0.
\end{aligned}
\end{equation}
Since this equation holds for any $m\in \pz$, the coefficients of all power of $m$ must be 0.
Next we prove the non-existence of $M$ with one type $A$ component in two cases.

{\em Case 1: $\beir=1$.}
If $t>2$, then the coefficient of $m^{t+2}$ in the right side of (\ref{eq4.35}) is
$$ra_t(k+1)^2(k^t+1-(k+1)^t)=-ra_t(k+1)^2\sum_{j=1}^{t-1}\binom tjk^j,$$
which is nonzero if $k\neq 0$ and $k\neq -1$, a contradiction.

If $t=2$ then (\ref{eq4.35}) turns to
$$\big(-r+(k+1)m\big)\big(k(k+1)(a_0-2ra_2)m^3-ra_0(k+1)^2m^2+ra_0(k+1)(1-a)m\big)=0,$$
which implies $a_0=a_2=0$, a contradiction.

If $t\leq 1$ then (\ref{eq4.35}) turns to
$$\big(-r+(k+1)m\big)\big(-ra_0(k+1)^2m^2+ra_0(k+1)(1-a)m\big)=0,$$
whose left hand side has coefficient of $m^4$
$$a_0 k(k+1)^2\neq 0\ \ \ \text{ if $k\neq 0$ and }k\neq -1.$$
With this contradiction we prove the case $\beir=1$ does not exist.

{\em Case 2: $\beir=0$.}
In this case equation (\ref{eq4.35}) turns to
\begin{equation}\label{eq4.36}
\begin{aligned}
-&r^2(k+1)m(a_0+a_1km+\cdots+a_tk^tm^t)-r^2(k+1)m(a_0+a_1m+\cdots+a_tm^t)\\
-&a_0r(k+1)m((k+1)m+a)(-r+km)+a_0rk(k+1)m^2(km+a)\\
+&r^2(k+1)m(a_0+a_1(k+1)m+\cdots+a_t(k+1)^tm^t)\\
+&r^2(-r+km)(a_0+a_1(k+1)m+\cdots+a_t(k+1)^tm^t)=0\text{ for } k\neq 0.
\end{aligned}
\end{equation}
If $t>2$ then the coefficient of $m^{t+1}$ in the left side of (\ref{eq4.36}) is
$$-r^2(k+1)a_tk^t-r^2(k+1)a_t+r^2 a_t(k+1)^t+r^2 a_t(k+1)^t
  =r^2(k+1)a_t\sum_{j=1}^{t-1}\binom tj k^j,$$
which is nonzero if $k\neq -1,0$. This is a contradiction.

If $t=2$ then the coefficient of $m^2$ in the left side of (\ref{eq4.36}) is
$r^2(k+1)^2a_0\neq0$ if $k\neq -1$. Also a contradiction.

If $t\leq 1$, then the coefficient of $m^3$ in the left side of (\ref{eq4.36}) is
$-rk(k+1)a_0\sum_{j=1}^{t-1}\binom tj k^j\neq 0$ if $k\neq -1,0$.
This contradiction disavows the case $\beir=0$.

In conclusion there does not exist module of intermediate series
for the algebra $\g$ with
more than one components among which one is of type $A$.

\subsection{All components being of type $V$ except one $B$}

In this subsection we deal with the case when $M$ has more than one components
of which one is of type $B$.
We prove such a module does not exist.

Assume the component $M(i)=B_i(a)$ for some $a\in\C$.
Then all other components are of type $V$.
Notice that the unanimous parameter $\al$ for all components equals to $-i$.
Let $i-s\rightsquigarrow i$.
We may write $M(i-s)=V_{i-s}(-i,\beims)$ for some $\beims\in\C$.

For $m,n\in \pz$ such that $m+n\neq 0$, consider
$$\begin{aligned}
&(n+s)\f {m+n+s}{i-s}\ve{i+m+n}=[L_m,L_{n+s}]\ve{i-s}\\
&=(n\f{n+s}{i-s}-(-s+\beims m)\f{n+s}{m+i-s})\ve{i+m+n}.
\end{aligned}$$
Hence
\begin{equation}\label{eq4.37}
(n+s)\f {m+n+s}{i-s}=n\f{n+s}{i-s}-(-s+\beims m)\f{n+s}{m+i-s}
    \text{ for } m+n\neq 0.
\end{equation}
Let $n=0$, we get
\begin{equation}\label{eq4.38}
-s\f{m+s}{i-s}=(-s+\beims m)\f{s}{m+i-s}\ \ \ \text{ for } m\neq 0,
\end{equation}
or
\begin{equation}\label{eq4.39}
s\f{-m+s}{i-s}=(s+\beims m)\f{s}{i-s-m}\ \ \ \text{ for } m\neq 0.
\end{equation}
Replace $m$ by $km$ and $n$ by $m$ in (\ref{eq4.37}), then we have
$$(-s+\beims km)\f{m+s}{km+i-s}=m\f{m+s}{i-s}-(m+s)\f {m+km+s}{i-s}
  \ \ \ \text{ for } m\neq 0, k\neq -1.$$

Consider $[L_m,L_s]\ve{i-s+km}$ for $k\neq -1$ and we have
\begin{equation}\label{eq4.40}
s\f{m+s}{i-s+km}=km\f {s}{i-s+km}-(-s+km+\beims m)\f{s}{i-s+m+km}
  \ \ \ \text{ for }k\neq -1.
\end{equation}
Now multiply $(-s+\beims km)(-s+\beims(k+1)m)$ to (\ref{eq4.40})
and apply (\ref{eq4.38}) and (\ref{eq4.39}), then we get
\begin{equation}\label{eq4.41}
\begin{aligned}
&\f{(k+1)m+s}{i-s}(-s+km+\beims m)(-s+\beims km)-\f{m+s}{i-s}m(-s+\beims(k+1)m)\\
&+\f{(k+1)m+s}{i-s}(m+s)(-s+\beims(k+1)m)-\f{km+s}{i-s}km(-s+\beims(k+1)m)=0\\
&\ \ \ \ \ \ \ \ \ \ \ \ \ \ \ \ \ \ \ \ \ \ \ \ \ \ \ \text{ for }k\neq -1.
\end{aligned}
\end{equation}
This implies that $\f{m+s}{i-s}$ is a polynomial in $m$.
Write
$$\f{m+s}{i-s}=\sum_{j=0}^{t}a_jm^j,$$
where $t=\partial\f{m+s}{i-s}$ is the order of $\f{m+s}{i-s}$,
$a_j\in\C, a_t\neq 0$ and $a_0=\f{s}{i-s}\neq 0$.

Consider $[L_m,L_s]\ve{i-s-m}$ and we get
\begin{equation}\label{eq4.42}
s\f{m+s}{i-s-m}=-m(m+a)\f {s}{i-s-m}+(s+m-\beims m)\f{s}{i-s}.
\end{equation}
Consider $[L_{-m},L_{m+s}]\ve{i-s}$ and we get
$$(s+\beims m)\f{m+s}{i-s-m}=m(m-a)\f{m+s}{i-s}+(m+s)\f s{i-s}.$$
Combine with (\ref{eq4.42}) and apply (\ref{eq4.39}) to give
\begin{equation}\label{eq4.43}
sm(m-a)\f{m+s}{i-s}+sm(m+a)\f{-m+s}{i-s}=\beims(1-\beims)m^2\f s{i-s}.
\end{equation}
Next we continue to discuss in two cases.

{\em Case 1: $\beims=0$.}
In this case (\ref{eq4.41}) turns to
$$(k+1)\f{(k+1)m+s}{i-s}-k\f{km+s}{i-s}-\f{m+s}{i-s}=0.$$
Then the leading coefficient in this equation is
$$a_t(k+1)^{t+1}-a_tk^{t+1}-a_t=\sum_{j=1}^t\binom{t+1}j k^j,$$
which is nonzero if $t>0$ and $k>0$.
This forces $t=0$ and $\f{m+s}{i-s}=\f{s}{i-s}$.
Then (\ref{eq4.43}) turns to
$$sm^2\f{s}{i-s}=0,$$
which is a contradiction.

{\em Case 2: $\beims\neq 0$.}
If $t>0$, then expending equation (\ref{eq4.41}) we get the leading coefficient
\begin{equation}\label{eq4.44}
-a_t\beims[k^{t+1}(\beims+t)+k^t(t\beims+\frac{t^2+t+2}{2})+\cdots+k(\beims+t)]=0,
\end{equation}
which stands for any $k\neq -1$.
If $t>1$, then we have
$$\beims+t=0,\ \ \ \ t\beims+\frac{t^2+t+2}{2}=0,$$
which forces $t=2$ and $\beims=-2$.
Then (\ref{eq4.43}) becomes
$$2sa_2m^4+2m^2(sa_0+3a_0-saa_1)=0,$$
which is a contradiction if we choose $m$ large enough.

If $t=1$ then by (\ref{eq4.44}) we have $\beims=-t=-1$.
Therefore (\ref{eq4.41}) turns to
$$2k(sa_1-a_0)m^2+sa_0(k-1)m=0\text{ for any $m$ and } k\neq -1.$$
Take $k=0$ and $m\neq 0$, then we get $a_0=0$,
which is a contradiction.
This forces $t=0$ and $\f{m+s}{i-s}=\f{s}{i-s}$.

Now (\ref{eq4.41}) becomes
$$k\beims(1-\beims)m^2=0,$$
which implies $\beims=1$.
Therefore, (\ref{eq4.43}) turns to
$$2sa_0m^2=0,$$
which is also a contradiction.

In conclusion there does not exist module of intermediate series
for the algebra $\g$ with
more than one components among which one is of type $B$.
Now we have completed the proof of Theorem \ref{thm4.3}.

\section{More about the matrix $F$ and examples of MOIS}
\def\theequation{5.\arabic{equation}}
\setcounter{equation}{0}

In this section we prove some more properties about
the MOIS over $\g$ and the corresponding matrix $F$.
At last, for small $p$, some explicit examples of matrix $F$ are given,
through which one should see more clear the structures of the matrix $F$
and a MOIS over $\g$.

Let $M$ be an irreducible MOIS over $\g$ with weight space decomposition
$M=\oplus_{i\in\Z}M_{i}$, where $M_i$ is the weight space
with respect to the weight $\al+i$ for some $\al\in\C$.
By Theorem \ref{thm4.3} we may write $M=\V$ for some $\be\in\C$.
Recall $\g''=\spanc{L_s\mid s\notin\pz}$ and the set $o(F)$.

\begin{prop}\label{prop5.1}
Assume that $M$ has more than one components.
Let $M(i)$ be a component of $M$ and $s\notin\pz$.\\
(1) $M=\U(\g'')M(i)$.\\
(2) If $L_sM(i)\neq 0$, then $L_sM=M$ and $L_sM(j)\neq 0$ for any component $M(j)$ of $M$.\\
(3) If $\f si\neq 0$ and $\mathrm{gcd}(s,p)=1$, then $o(F)=\{0,1,2,\cdots, p-1\}$.\\
(4) If $p$ is a prime, then $o(F)=\{0,1,2,\cdots, p-1\}$.\\
(5) If $\f si\neq0, \f rj\neq0$ and $\mathrm{gcd}(r,s)=1$, then $o(F)=\{0,1,2,\cdots, p-1\}$.
\end{prop}
\begin{pf}
Notice that $\U(\g'')M(i)$ is a nonzero $\g$-submodule of $M$.
Then (1) is clear since $M$ is irreducible.
For (2), notice that central elements of $\g$ acts on $M$ trivially.
Hence $L_sM$ is a nonzero $\g$-submodule of $M$.
Therefore, $L_sM=M$.
Suppose $L_sM(j)=0$ for some component $M(j)$ of $M$. Then
$$L_sM=L_s\U(\g'')M(j)=\U(\g'')L_s M(j)=0,$$
which is a contradiction.
(3), (4) and (5) follow from (2) and the condition (II).
\pfend

Now we list some examples.
All $\f ij$ appearing in the following examples are assumed to be nonzero complex numbers.

\noindent
{\bf Example 1}: $p=3$. Set $F=\begin{pmatrix}
\f 10 &\f 11 &\f 12\\
0     &  0   &  0
\end{pmatrix}$ and
$\V=V(0)\oplus V(1)\oplus V(2)$.
Here $L_1$ links all $V(i)$, but $L_2$ acts on $\V$ trivially.
The linkage in $\V$ is
$$\xymatrix@!0{
 0 \ar[rr] &          & 1 \ar[dl] \\
           & 2\ar[ul] &
}.$$
For any nonzero $\f 1j, j=0,1,2$,
since the condition (III) stands,
a matrix of such form may be equipped with a MOIS over $\g$.

\noindent
{\bf Example 2}: $p=4$. Set $F=\begin{pmatrix}
  0   &  0  &  0    &  0\\
\f 20 &  0  & \f 22 &  0\\
0     &  0   &  0   &  0
\end{pmatrix}$,
which satisfies the conditions (I)-(III).
Then there exists a MOIS $\V=V(0)\oplus V(2)$ over $\g$.
Here $L_2$ links $V(0)$ and $V(2)$ to each other, but $L_1,L_3$ act on $\V$ trivially.
This example shows the existence of a MOIS whose number of components is less than $p$.

\noindent
{\bf Example 3}: $p=5$. Set $F=\begin{pmatrix}
\f 10 &\f 11 & \f 12 &\f 13 &\f 14\\
\f 20 &\f 21 & \f 22 &\f 23 &\f 24\\
0     &  0   &  0    &  0   &  0  \\
0     &  0   &  0    &  0   &  0
\end{pmatrix}$.
In order for $F$ to equip with a MOIS over $\g$,
$F$ has to satisfy the condition (III),
which forces
\begin{equation}\label{eq5.1}
 \f21=\frac{\f12\f20}{\f10};\ \ \f22=\frac{\f12\f13\f20}{\f10\f11};\ \
  \f23=\frac{\f13\f14\f20}{\f10\f11};\ \ \f24=\frac{\f14\f20}{\f11}.
 \end{equation}
Now for a matrix $F$ satisfying (\ref{eq5.1}) with arbitrary nonzero $\f20,\f1i,0\leq i\leq 4$,
we can equip with it a MOIS $\V=\bigoplus_{i=0}^4 V(i)$ over $\g$.
Here $L_1, L_2$ link all $V(i)$, while $L_3,L_4$ act on $\V$ trivially.
The linkage here is
$$
\xymatrix@!0{
   &    &  0  &     \\
   &  2 &     &  3  \\
   &    &     &     \\
 1 &    &  4  &     \\
\ar "1,3"; "4,1"   
\ar "4,1"; "2,2"
\ar "2,2"; "2,4"
\ar "2,4"; "4,3"
\ar "4,3"; "1,3"
\ar "1,3"; "2,2"
\ar "2,2"; "4,3"
\ar "4,3"; "4,1"
\ar "4,1"; "2,4"
\ar "2,4"; "1,3"
}\ .$$

\noindent {\bf Example 4}: Let $p=4$ and
$F=\begin{pmatrix}
\f 10 & \f11 & \f12  &\f13\\
\f 20 &  0   & \f 22 &  0 \\
  0   &  0   &  0    &  0
\end{pmatrix}$.
Such a matrix $F$ can not be equipped with a MOIS.
Suppose, otherwise, that $F$ is equipped with a MOIS $M$.
Since $M(1)\neq0$ is a component of $M$ by $\f11\neq 0$,
we have $\f21\neq 0$ by $\f20\neq0$ and Proposition \ref{prop5.1}(3),
which is a contradiction.

\noindent
{\bf Example 5}: Let $p=12$ and $F$ be such that all entries are zero except $\f 80,\f84,\f88$.
Such $F$ satisfies the conditions (I)-(III) and
the corresponding MOIS is $M=M(0)\oplus M(4)\oplus M(8)$ with linkage
$$
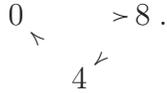
\xymatrix@!0{
 0 \ar[rr] &          & 8 \ar[dl] \\
           & 4\ar[ul] &
}.$$
This example shows that even if $o(F)$ is a proper subset of $\{0,1,2,\cdots, p-1\}$
and $\f si\neq0$, we may not have $s\mid p$.


\begin{thebibliography}{99}


\bibitem[ACKP]{ACKP}
E. Arbarello, C. De Concini, V.G. Kac, C. Procesi,
Moduli spaces of curves and representation theory,
Comm. Math. Phys. 117(1988) 1-36.


\bibitem[B]{B}
Y. Billig, Jet modules,
Canad. J. Math. Vol. 59 (4) (2007) 712-729.



\bibitem[BF]{BF}
Y. Billig, V. Futorny,
Classification of irreducible representations of Lie algebra of vector fields on a torus,
J. Reine Angew. Math. 720 (2016) 199-216.


\bibitem[BL]{BL}
D. Britten, F. Lemire,
On level 0 affine Lie modules,
Canad. Math. Bull. 37 (3) (1994) 310-314.


\bibitem[BGK]{BGK}
S. Berman, Y. Gao, Y. S. Krylyuk,
Quantum tori and elliptic quasi-simple Lie algebras,
J. Funct. Anal. 135 (2) (1996) 339-389.


%


\bibitem[GGS]{GGS}
M. Gao, Y. Gao, Y. Su,
Irreducible quasi-finite representations of a Block type Lie algebra,
Commun. Alg. 42 (2014) 511-527.



\bibitem[LZ1]{LZ1}
R. Lu, K. Zhao,
Classification of irreducible weight modules over higher rank Virasoro algebras,
Adv. Math. 206 (2006) 630-656.


\bibitem[LZ2]{LZ2}
R. Lu, K. Zhao,
Classification of irreducible weight modules over the twisted Heisenberg-Virasoro algebra,
Commun. Contemp. Math. 12 (2) (2010) 183-205.



\bibitem[M]{M}
O. Mathieu,
Classification of Harish-Chandra modules over the Virasoro Lie algebra,
Invent. Math. 107 (1992) 225-234.


\bibitem[MP]{MP}
C. Martin, A. Piard,
Indecomposable modules over the Virasoro Lie algebra and a conjecture of V. Kac,
Commun. Math. Phys. 137 (1991) 109-132.


\bibitem[SZ]{SZ}
Y. Su, K. Zhao,
Generalized Virasoro and super-Virasoro algebras and modules of intermediate series,
J. Algebra 252 (2002) 1-19.


\bibitem[WT]{WT}
Q. Wang, S. Tan,
Quasifinite modules of a Lie algebra related to Block type,
J. Pure Appl. Alg. 211 (2007) 596-608.

\end{thebibliography}
\end{document}